\documentclass[10pt]{article}

\usepackage[T1]{fontenc}
\usepackage{lmodern}
\usepackage[utf8x]{inputenx}
\usepackage{microtype}
\usepackage{framed}
\usepackage{listings}
\usepackage{vmargin}
\usepackage{setspace}
\usepackage{mathrsfs, mathenv}
\usepackage{amsmath, amsthm, amssymb, amsfonts, amscd}
\usepackage{graphicx}
\usepackage{epstopdf}
\usepackage[svgnames]{xcolor}
\usepackage{hyperref}
\usepackage[capitalise]{cleveref}
\usepackage{xparse}
\hypersetup{citecolor=blue, colorlinks=true, linkcolor=black}
\usepackage{subcaption}
\usepackage{bbm}
\usepackage{cite}
\usepackage{verbatim}
\usepackage{pgfplots}
\usepackage{tikz}
\usepackage{color}
\usepackage{stmaryrd}
\usepackage{autonum}

\theoremstyle{plain}
\newtheorem{theorem}{Theorem}
\newtheorem{corollary}{Corollary}
\newtheorem{lemma}{Lemma}

\theoremstyle{definition}
\newtheorem{definition}{Definition}

\theoremstyle{remark}
\newtheorem{remark}{Remark}
\newtheorem{assumption}{Assumption}
\newtheorem{example}{Example}

\ifpdf
  \DeclareGraphicsExtensions{.eps,.pdf,.png,.jpg}
\else
  \DeclareGraphicsExtensions{.eps}
\fi

\usepackage{mathtools}

\usepackage{booktabs}

\usepackage{pgfplots}
\usepackage{tikz}
\usetikzlibrary{arrows,decorations.pathmorphing,backgrounds,positioning,fit,matrix}
\usepackage[labelfont=bf]{caption}
\usepackage{here}
\usepackage[font=normal]{subcaption}

\usepackage{enumitem}
\setlist[itemize]{leftmargin=.5in}
\setlist[enumerate]{leftmargin=.5in,topsep=3pt,itemsep=3pt,label=(\roman*)}

\newcommand{\email}[1]{\href{#1}{#1}}
\newcommand{\TheTitle}{Random time step probabilistic methods for uncertainty quantification in chaotic and geometric numerical integration} 
\newcommand{\TheAuthors}{A. Abdulle, G. Garegnani}
\title{{\TheTitle}}
\author{Assyr Abdulle\thanks{Mathematics Section, \'Ecole Polytechnique F\'ed\'erale de Lausanne (\email{assyr.abdulle@epfl.ch})}
	\and
	Giacomo Garegnani\thanks{Mathematics Section, \'Ecole Polytechnique F\'ed\'erale de Lausanne (\email{giacomo.garegnani@epfl.ch})}}
\date{}

\usepackage{amsopn}

\DeclarePairedDelimiter{\abs}{\lvert}{\rvert}
\DeclarePairedDelimiter{\norm}{\|}{\|}
\renewcommand{\phi}{\varphi}
\renewcommand{\theta}{\vartheta}

\renewcommand{\tilde}{\widetilde}
\renewcommand{\hat}{\widehat}

\newcommand{\iid}{\ensuremath{\stackrel{\text{i.i.d.}}{\sim}}}

\newcommand{\rightarrowtext}[1]{\ensuremath{\stackrel{#1}{\longrightarrow}}}
\newcommand{\leftrightarrowtext}[1]{\ensuremath{\stackrel{#1}{\longleftrightarrow}}}
\newcommand{\pdv}[2]{\ensuremath\partial_{#2}#1}
\newcommand{\N}{\mathbb{N}}
\newcommand{\R}{\mathbb{R}}

\newcommand{\OO}{\mathcal{O}}
\newcommand{\epl}{\varepsilon}
\newcommand{\diffL}{\mathcal{L}}

\newcommand{\defeq}{\coloneqq}

\newcommand{\Var}{\operatorname{Var}}
\newcommand{\E}{\operatorname{\mathbb{E}}}
\newcommand{\MSE}{\operatorname{MSE}}
\newcommand{\trace}{\operatorname{tr}}

\newcommand{\Hell}{d_{\mathrm{Hell}}}
\newcommand{\sksum}{{\textstyle\sum}}
\newcommand{\dd}{\,\mathrm{d}}
\definecolor{shade}{RGB}{100, 100, 100}
\definecolor{bordeaux}{RGB}{128, 0, 50}

\definecolor{leg1}{RGB}{0,114,189}
\definecolor{leg2}{RGB}{217,83,25}
\definecolor{leg3}{RGB}{237,177,32}
\definecolor{leg4}{RGB}{126,47,142}
\definecolor{leg5}{RGB}{119,172,48}

\definecolor{leg21}{RGB}{62,38,169}
\definecolor{leg22}{RGB}{46,135,247}
\definecolor{leg23}{RGB}{55,200,151}
\definecolor{leg24}{RGB}{254,195,56}

\ifpdf
\hypersetup{
	pdftitle={\TheTitle},
	pdfauthor={\TheAuthors}
}
\fi

\begin{document}
\maketitle	

\begin{abstract} 
A novel probabilistic numerical method for  quantifying the uncertainty induced by the time integration of ordinary differential equations (ODEs) is introduced. Departing from the classical strategy to randomize ODE solvers by adding a random forcing term, we show that a probability measure over the numerical solution of ODEs can be obtained by introducing suitable random time-steps in a classical time integrator.
This  intrinsic randomization allows for the conservation of geometric properties of the underlying deterministic integrator such as mass conservation, symplecticity or conservation of first integrals. Weak and mean-square convergence analysis are derived. We also analyse the convergence of the Monte Carlo estimator for the proposed random time step method and show that the measure obtained with repeated sampling converges in the mean-square sense independently of the number of samples. Numerical examples including chaotic Hamiltonian systems, chemical reactions and Bayesian inferential problems illustrate the accuracy, robustness and versatility of our probabilistic numerical method.
\end{abstract}



\normalsize
\section{Introduction} 
A variety of methods for integrating ordinary differential equations (ODEs) has been studied in the last decades, \cite{HNW93, HaW96, HLW06}, with an emphasis on building accurate and stable deterministic approximations of the exact solution. In general, these methods are based on a time discretization on which the solution of the ODE is approximated via an iterative deterministic algorithm. Given a time step $h$, which indicates the refinement of the discretization, all these methods provide a point value for the approximation of the solution and guarantee that in the asymptotic limit of $h \to 0$ the numerical approximation will coincide with the exact solution. However, for some problems such as chaotic systems or inference problems having a distributional solution can help to quantify the uncertainty introduced by the numerical discretization without invoking the asymptotic limit $h \to 0$.

In recent years, probabilistic numerical methods for differential equations have been proposed \cite{CGS17, CCC16, SDH14} in order to quantify the uncertainty introduced by the time discretization in a statistical manner. A review summarizing the recent advancements in the field of probabilistic numerical can be found in \cite{OaS19, COS19}. In general, these methods proceed iteratively to establish a probability measure over the numerical solution, thus providing a richer information than a single point value. In particular, probabilistic solvers offer a quantitative characterisation of late time errors by tuning the noise introduced by the method according to the accuracy of the solver. In this way, it is possible to obtain a reliable approach for capturing the sensitivity of the solution to numerical error, while transferring the convergence properties of classical deterministic integrators to the introduced probability measure in a consistent manner. 

In the following, we will first show two examples motivating the probabilistic approach, and then present the main contributions of this work.

\subsection{Motivating examples} Probabilistic integrators for ODEs do not provide more accurate solutions than classical deterministic methods nor are they computationally cheaper. Nevertheless, they can be useful in a variety of different problems, among which we identified the integration of chaotic dynamical systems and the solution of Bayesian inverse problems, which are briefly presented here.

\subsubsection*{Chaotic differential equations} 
\begin{figure}
	\begin{center}
		\begin{tabular}{c}
			\includegraphics[]{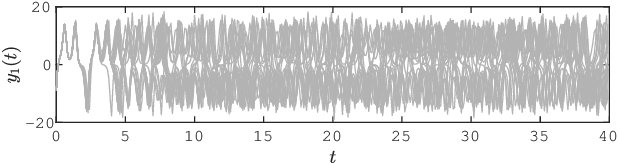} \\
			\includegraphics[]{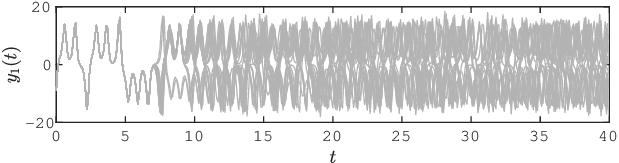} \\
			\includegraphics[]{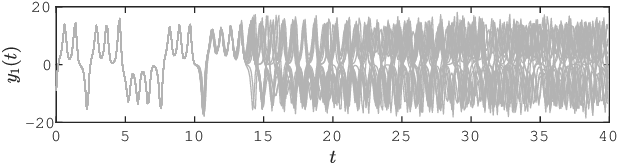}
		\end{tabular}
	\end{center}
	\caption{First component {$y_1(t)$} of the solution of \eqref{eq:Lorenz} with decreasing Gaussian perturbations on the initial condition from top to bottom {($\sigma = 10^{-1}, 10^{-3}, 10^{-5}$, respectively)}.}
	\label{fig:LorenzTest}
\end{figure}

Let us consider the Lorenz system \cite{Lor63}, which is defined by the following ODE
\begin{equation}\label{eq:Lorenz}
\begin{aligned}
	y_1' &= \eta(y_2 - y_1), \quad &&y_1(0) = -10,\\
	y_2' &= y_1(\rho - y_3) - y_2, \quad &&y_2(0) = -1,\\
	y_3' &= y_1y_2 - \beta y_3, \quad &&y_3(0) = 40.
\end{aligned}
\end{equation}
It is well-known that for $\rho=28$, $\eta=10$, $\beta=8/3$, this equation has a chaotic behaviour, i.e., a small perturbation forces the trajectories to deviate from the true solution. Integrating numerically \eqref{eq:Lorenz} the error which is introduced at each time step is indeed a perturbation, thus any numerical solution cannot be considered reliable. In order to explore the state space of this chaotic dynamical system, we introduce a random perturbation on the initial condition, implemented as a scalar Gaussian random variable $\epl \sim \mathcal{N}(0, \sigma^2)$ and artificially added to the first component ${y_1(t)}$ at time $t = 0$. In \cref{fig:LorenzTest} we show $M = 20$ numerical trajectories given by a second-order Runge--Kutta method for three different scales of the noise. It is possible to remark that in all the three cases, the numerical solutions almost coincide up to some time $\bar t$, thus diverging and showing the chaotic nature of the Lorenz system. It could be argued that up to time $\bar t$, the numerical solution offers a reliable approximation of the true solution as the dynamics have not yet switched to the chaotic regime. Nevertheless, it is unclear how to choose ${\sigma^2}$ so that the amount of noise that is introduced is balanced with the numerical error. Probabilistic methods for differential equations such as the one presented in this work and the one introduced by Conrad et al. \cite{CGS17} provide a rigorous analysis that suggests how to introduce a source of artificial noise in a consistent manner.
 
\subsubsection*{Bayesian inference} Problems of Bayesian inference are most often used to justify the usefulness of probabilistic methods for differential equations. The impact of a probabilistic component in the numerical approximation of inverse problems involving ODEs has already been presented in several works (e.g., \cite{CGS17, CCC16, COS17}). In particular, the common underlying idea of these works is that if a deterministic integrator with a fixed finite time step is employed to approximate the solution of the ODE appearing in an inverse problem, the numerical error introduced by deterministic solvers can lead to inappropriate and non-predictive posterior concentrations. In the limit of an infinitely refined time discretization the posterior distributions obtained with a classical numerical method will indeed tend to the true distribution, but for a fixed time step (i.e., for a fixed computational budget) numerical error can lead to posterior concentrations away from the true value of the parameter of interest. These inappropriate solutions to inverse problems can be corrected by employing a probabilistic integrator to solve the ODE, thus obtaining posterior distributions that reflect the uncertainty given by the numerical solver (see \cref{fig:AnalyticalPosterior} at page \pageref{fig:AnalyticalPosterior} for an example).

\subsection{Contributions} The method we analyse in this paper is inspired from the work of Conrad et al. \cite{CGS17}, where a probabilistic method for ODEs is presented. This method consists of perturbing a deterministic numerical solution (e.g. arising from a Runge--Kutta discretization) with an additive source of noise at each time step. By appropriately scaling the random term, they manage to obtain a probabilistic solution without altering the convergence of the underlying deterministic scheme. 

An additive noise contribution could nonetheless produce disruptive effects on favourable geometric features of deterministic schemes. A direct example of this non-robust behaviour is given by ODEs for which the solution is supposed to stay positive and small. In this case, the addition of a random contribution could force the solution in the negative plane, hence the numerical solution could be physically meaningless. Chemical reactions with small population size for one species at some time of the evolution are typical physical examples. In particular, an additive random term could force the solution on the negative plane with a non-zero probability, and this probability could become non-negligibly big in case the magnitude of one component of the solution is small. Other geometric properties of an underlying ODE are also destroyed when perturbing the flow by a noisy forcing term.

Motivated by these issues, we present in this work a new probabilistic method for ODEs based on a random selection of the time steps. Hence, the randomness of the scheme becomes intrinsic in contrast to the additive noise method. For this new robust probabilistic integrator, we are able to prove strong and weak convergence towards the exact solution of the underlying ODE. Precisely, setting the variance of the random time steps to be proportional to some power of a deterministic time step allows to retrieve the rates of the underlying Runge--Kutta integrator.

It has been pointed out by Kersting and Hennig \cite{KeH16} that probabilistic methods based on sampling should be equipped with a criterion to choose the number of samples, so that computational effort is not wasted or, conversely, the sample size is not insufficient to describe the dynamics in a probabilistic fashion. In order to address this issue, in this work we show that Monte Carlo estimators drawn from our probabilistic solver converge with respect to the time step in the mean square sense independently of the sample size. We are able to prove a similar property for the scheme proposed in \cite{CGS17}.

A large variety of dynamical systems is characterised by geometrical properties of their flow map \cite{HLW06}. Most notably, Hamiltonian systems, which are employed for modelling a variety of physical phenomena, are endowed with the property of symplecticity. It is possible to obtain good approximations of the solutions of Hamiltonian systems via mimicking numerically the geometric properties of the exact flow, i.e., by employing symplectic integrators. In particular, for symplectic integrators the energy function conserved by the exact flow is approximately conserved by numerical trajectories over long time spans, which in turn guarantees high-quality numerical solutions at the price of a rather low computational effort. While geometric properties of Runge--Kutta schemes have been analysed extensively in the deterministic case, they have not been considered yet for probabilistic numerical methods. The method we present in this work, being only an intrinsic modification of a Runge--Kutta integrator, is endowed with the geometric properties of its deterministic counterpart. In particular, we first show that our probabilistic scheme inherits the property of exact conservation of first integrals of the considered dynamics. Then, we show that in Hamiltonian systems the good approximation of the energy function given by symplectic schemes is preserved by our randomisation procedure over polynomially long times.

\subsection{Outline} The paper is organised as follows. In \cref{sec:MethodIntro} we introduce the setting for probabilistic numerics and present our novel numerical scheme. We then show in \cref{sec:WeakOrder} and \cref{sec:StrongOrder} the properties of weak and mean square convergence of the numerical solution towards the exact solution of the ODE. In \cref{sec:MonteCarlo} we analyse the accuracy of Monte Carlo estimators drawn from the numerical solution. The geometric properties of the numerical scheme are presented in \cref{sec:GeomProperties} and \cref{sec:Hamiltonian}, while in \cref{sec:BayesianInference} we introduce Bayesian inverse problems in the ODE setting, and show how our method can be integrated in existing sampling strategies. Finally, we show a variety of numerical experiments confirming our theoretical results in \cref{sec:NumericalExperiments}.

\section{Random time step Runge--Kutta method}\label{sec:MethodIntro}

\begin{figure}
\centering
\begin{tikzpicture}
	\draw[->] (-0.3,0) -- (8,0) node[right] {\small $t$};
	\draw[->] (0,-0.3) -- (0,3.5);
	\draw[scale=1,domain=0.5:4.2,smooth,variable=\t,blue] plot ({\t},{0.2 + 0.15*\t*\t}) node[right] {\small $y(t)$};
	\draw[scale=1,domain=1:6,smooth,variable=\t,black] plot ({\t},{0.05 + 0.3*\t});
	\draw (1,0.35) node[circle,draw=black, fill=white, inner sep=0pt,minimum size=3pt]{};
	\draw[dashed] (1,0) -- (1,0.35) node[above]{\small $y_0$};
	\draw[dashed] (4,0) -- (4,1.25) node[below right] {\small $y_1$};
	\draw (4,1.25) node[circle,draw=black, fill=white, inner sep=0pt,minimum size=3pt]{};
	\draw[-latex] (6,1.85) -- (4.05,1.85) node[left] {\small $Y_1$};
	\draw (4,1.85) node[circle,draw=black, fill=white, inner sep=0pt,minimum size=3pt]{};
	\draw[dashed] (6,0) -- (6,1.85);
	\draw (6,1.85) node[circle,draw=black, fill=white, inner sep=0pt,minimum size=3pt]{};
	\draw (1,0) node[below] {\small $t_0$};
	\draw (4,0) node[below] {\small $t_0 + h$};
	\draw (6,0) node[below] {\small $t_0 + H_0$};
	\draw[latex'-latex',red] (4,1.3) -- (4,1.8);
\end{tikzpicture}
\caption{Graphical representation of one step of the RTS-RK method with $\Psi_h(y) = y + hf(y)$. The red arrow is the stochastic contribution due to random time-stepping.}
\label{fig:GraphRandomStep}
\end{figure}
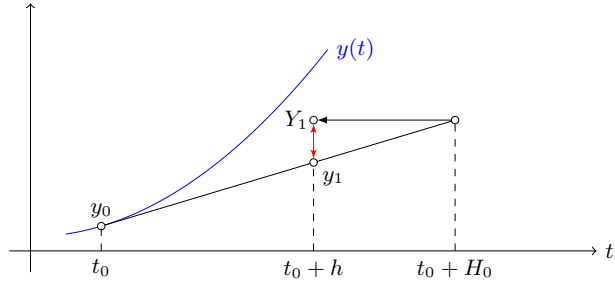

Let us consider a Lipschitz function $f\colon\R^d\to\R^d$ and the ODE
\begin{equation}\label{eq:ODE}
	y' = f(y), \quad y(0) = y_0 \in \R^d.
\end{equation}
In the following, we will write for simplicity the solution $y(t)$ of \eqref{eq:ODE} in terms of the flow of the ODE. In particular, we consider the family  $\{\phi_t\}_{t \geq 0}$ of functions $\phi_t\colon\R^d\to\R^d$ such that 
\begin{equation}
	y(t) = \phi_t(y_0).
\end{equation}
Given a time step $h$, let us consider a Runge--Kutta method which deterministically approximates the solution $\phi_t(y_0)$ of \eqref{eq:ODE}. In particular, we can write the numerical solution $y_k$ approximating $\phi_{t_k}(y_0)$, with $t_k = kh$ in terms of the numerical flow $\{\Psi_t\}_{t \geq 0}$, with $\Psi_t\colon\R^d\to\R^d$, which is uniquely determined by the coefficients of the method, as
\begin{equation}
	y_{k+1} = \Psi_h(y_k), \quad k = 0, 1, \ldots.
\end{equation}
In order to provide a probabilistic interpretation of the numerical solution rather than a series of point values, Conrad et al. propose the scheme defined by
\begin{equation}\label{eq:AN-RK}
	Y_{k+1} = \Psi_h(Y_{k}) + \xi_k(h), \quad k = 0, 1, \ldots,
\end{equation}
where $Y_k$ is a random variable approximating $y(t_k)$ with $Y_0 = y_0$, and $\xi_k(h)$ are appropriately scaled independent and identically distributed (i.i.d.) random variables with values in $\R^d$. Maintaining the same notation as in \eqref{eq:AN-RK}, in this work we propose a random time-stepping Runge--Kutta method (RTS-RK), i.e., the scheme defined defined by the recurrence relation
\begin{equation}\label{eq:RTS-RK}
	Y_{k+1} = \Psi_{H_k}(Y_k), \quad k = 0, 1, \ldots,
\end{equation}
where $Y_k$ is still a random variable approximating $y(t_k)$ and the time steps $H_k$ are locally given by a sequence  of i.i.d. random variables with values in $\R^+$. A graphical representation of one step of the RTS-RK method is given in \cref{fig:GraphRandomStep}. Let us finally remark that the sequence $Y_k$, $k = 0, 1, \ldots$, form a homogeneous Markov chain, as the transition probability is independent of the index $k$.
\begin{remark} We note that in terms of computational cost simulating the two methods \eqref{eq:AN-RK} and \eqref{eq:RTS-RK} is equivalent.
\end{remark}

\subsection{Assumptions and notation} We now present the main assumptions and notations which are used throughout the rest of this work. Firstly, we have to consider the possible values taken by the random step sizes, which have to satisfy restrictions that are necessary not to spoil the properties of deterministic methods. 
\begin{assumption}\label{as:hStrong} The i.i.d. random variables $H_k$ satisfy for all $k = 0, 1, \ldots$
	\begin{enumerate}
		\item\label{as:hStrong_Pos} $H_k > 0$ a.s.,
		\item\label{as:hStrong_E} there exists $h > 0$ such that $\E H_k = h$,
		\item\label{as:hStrong_Var} there exist $p \geq 1/2$ and $C > 0$ independent of $k$ such that the scaled random variables $Z_k \defeq H_k - h$ satisfy
		\begin{equation}
			\E Z_k^2 = Ch^{2p+1}.
		\end{equation}
	\end{enumerate}
\end{assumption}
The class of random variables satisfying the hypotheses above is general. However, it is practical for an implementation point of view to have examples of these variables.
\begin{example}\label{ex:uniformH} Let us consider the random variables $\{H_k\}_{k\geq 0}$ such that
	\begin{equation}
	H_k \iid \mathcal{U}(h-h^{p+1/2}, h+h^{p+1/2}), \quad 0 < h < 1, \quad p \geq 1/2.
	\end{equation}
	We easily verify that the assumptions \ref{as:hStrong_Pos} and \ref{as:hStrong_E} are verified as $h < 1$, and that \ref{as:hStrong_Var} is verified with $C = 1/3$. Another choice of random variables could simply be 
	\begin{equation}\label{eq:logNormalH}
	H_k \iid \log\mathcal{N}\big(\log h - \log\sqrt{1 + h^{2p}}, \log(1 + h^{2p})\big),
	\end{equation}
	for which the properties above are trivially verified (with $C = 1$), provided $p > 1/2$.
\end{example}
We secondly introduce an assumption on the deterministic  method underlying the RTS-RK scheme, identified by its numerical flow $\Psi_h$. 
\begin{assumption}\label{as:PsiStrong} The Runge--Kutta method defined by the numerical flow $\{\Psi_t\}_{t\geq 0}$ is of order $q$, i.e., for $h$ small enough, there exists a constant $C > 0$ such that
		\begin{equation}
		\norm{\Psi_h(y) - \phi_h(y)} \leq Ch^{q+1}, \quad \forall y \in \R^d.
		\end{equation}
\end{assumption}

\begin{remark} Depending on the domain of definition of the vector field $f$, the choice of an unbounded distribution for the time step could give rise to two critical issues. In particular, 
	\begin{enumerate}
		\item if $f\colon D \to \R^d$, where $D \in \R^d$ is a bounded open subset of $\R^d$, allowing the time step to assume unbounded values as, e.g., in case of the log-normal distribution \eqref{eq:logNormalH}, may force the solution outside $D$,
		\item if $\Psi_h$ is the numerical flow of an implicit method, the solution could be ill-posed.
	\end{enumerate}
	In both the two cases above, we suggest to employ uniform time steps as in \cref{ex:uniformH}, which allow the time steps to be small enough almost surely. For the first issue, more sophisticated techniques of path rejection could be employed \cite{MiT05}, but the mean-square convergence properties which will be examined in \cref{sec:StrongOrder} would not hold. 
\end{remark}

In order to tackle the second issue presented in the Remark above, we introduce a further assumption.

\begin{assumption}\label{as:hImplicit} If the map $\Psi_t$ is implicit, the time steps $H_k$ satisfy $H_k \leq M < \infty$ almost surely, where $M$ is small enough to allow the scheme to be well-posed.
\end{assumption}

Let us finally remark that the choice of the distribution of the time steps is artificial and therefore arbitrary. Hence, choosing a bounded distribution does not represent a limitation to the numerical scheme.

\section{Weak convergence analysis}\label{sec:WeakOrder}

The first property of the RTS-RK method we wish to analyze is its weak convergence, which gives an indication about the behavior of the numerical solution \eqref{eq:RTS-RK} in the mean sense. In the following, we denote by $\mathcal C^l_b(\R^d, \R)$ the functions in $\mathcal C^l(\R^d, \R)$ with all derivatives up to order $l$ bounded uniformly in $\R^d$. Moreover, we consider the integration of \eqref{eq:ODE} over the finite length domain $[0,T]$, where $T > 0$ is the final time. Let us define the weak order of convergence. 
\begin{definition} The numerical method \eqref{eq:RTS-RK} has weak order $r$ for \eqref{eq:ODE} if for any sufficiently smooth function $\Phi\colon\R^d \to \R$ there exists a constant $C > 0$ independent of $h$ such that
	\begin{equation}
		\abs{\E\Phi(Y_k) - \Phi(y(kh))} \leq Ch^r,
	\end{equation}
	for all $k = 1, 2, \ldots, N$ and $T = Nh$.
\end{definition} 
Let us introduce the Lie derivative of the flow $\diffL = f \cdot \nabla$, which allows us to adopt the semi-group notation for the exact solution of \eqref{eq:ODE} (see e.g. \cite[Section III.5.1]{HLW06} or \cite[Section 4.3]{PaS08}) and write for any smooth function $\Phi$
\begin{equation}\label{eq:LieDerivative}
	\Phi(\phi_h(y)) = e^{h\diffL} \Phi(y).
\end{equation}
Moreover, let us recall that the probabilistic numerical solution $\{Y_k\}_{k\geq 0}$ forms a homogeneous Markov chain. Therefore, given $h > 0$ there exists an operator $\mathcal{P}_h$, the generator \cite[Section 2.3]{Pav14}, such that
\begin{equation}
	\E \big(\Phi(Y_{k+1})\mid Y_k = y\big) = (\mathcal{P}_h\Phi)(y).
\end{equation}
In order to have an analogy with the notation \eqref{eq:LieDerivative}, we adopt the exponential form of the infinitesimal generator and denote in the following $\mathcal{P}_h = e^{h\diffL_h}$, where we explicitly write the dependence of the Markov generator on the step size $h$. Furthermore, due to the homogeneity of the Markov chain, we can write
\begin{equation}\label{eq:MarkovProperty}
	\E \big(\Phi(Y_{k+1})\mid Y_0 = y\big) = e^{h\diffL_h}\E \big(\Phi(Y_k)\mid Y_0 = y\big).
\end{equation}
We can now state a result of local weak convergence of the probabilistic numerical solution.
\begin{lemma}[Weak local order]\label{lem:WeakLocalOrder} Let \cref{as:hStrong}, \cref{as:PsiStrong} and \cref{as:hImplicit} hold and let $f$ in \eqref{eq:ODE} be sufficiently smooth. If $\E|H_0^4| < \infty$, there exists a constant $C > 0$ independent of $h$ and $y$ such that for any function $\Phi\in \mathcal C^l_b(\R^d, \R)$, with $l = \max\{q, 3\}$
	\begin{equation}       
		\abs{\E(\Phi(Y_1)\mid Y_0 = y) - \Phi(\phi_h(y))} \leq C h^{\min\{2p+1, q+1\}}.
	\end{equation}
\end{lemma}
\begin{proof} Since $f$ is sufficiently smooth, the map $t \mapsto \Psi_t(y)$ is of class $C^2(\R^+, \R^d)$ and Lipschitz continuous with constant $L_\Psi$ independent of $y$. Let us expand the functional $\Phi$ computed on the numerical solution as
	\begin{equation}\label{eq:TaylorNumericalSolution}
		\begin{aligned}
			\Phi(Y_1) &= \Phi(\Psi_{H_0}(Y_0)) \\
			&= \Phi\Big(\Psi_h(Y_0) + (H_0-h)\partial_t\Psi_h(Y_0) + \frac{1}{2}(H_0-h)^2\partial_{tt}\Psi_h(Y_0) + \OO(\abs{H_0 - h}^3)\Big)\\
			&= \Phi(\Psi_h(Y_0)) + \Big((H_0 - h)\partial_t\Psi_h(Y_0)+\frac{1}{2}(H_0-h)^2\partial_{tt}\Psi_h(Y_0)\Big) \cdot \nabla\Phi(\Psi_h(Y_0))\\
			&\quad + \frac{1}{2}(H_0 - h)^2 \partial_t \Psi_h(Y_0) \partial_t \Psi_h(Y_0)^\top : \nabla^2\Phi(\Psi_h(Y_0)) + \OO(\abs{H_0 - h}^3),
		\end{aligned}
	\end{equation}
	where we denote by $\nabla^2\Phi$ the Hessian matrix of $\Phi$, and by $:$ the inner product on matrices induced by the Frobenius norm on $\R^d$, i.e., $A : B = \trace(A^\top B)$. Taking the conditional expectation with respect to $Y_0 = y$ and applying \cref{as:hStrong} we get
	\begin{equation}\label{eq:TaylorNumericalSolution2}
		\begin{aligned}
			e^{h\diffL_h}\Phi(y) - \Phi(\Psi_h(y)) &= \frac{1}{2} Ch^{2p+1}\partial_{tt}\Psi_h(y)\cdot \nabla\Phi(\Psi_h(y))\\
			&\quad + \frac{1}{2} Ch^{2p+1}\partial_t \Psi_h(y) \partial_t \Psi_h(y)^\top  \colon \nabla^2\Phi(\Psi_h(y)) + \OO(h^{3p+3/2}),
		\end{aligned}
	\end{equation}
	where we exploited Hölder's inequality for the last term. Moreover, expanding $\Phi$ in $y$ we get
	\begin{equation}
		\begin{aligned}
			\Phi(\Psi_h(y)) &= \Phi\left(\Psi_0(y) + h\partial_t \Psi_0(y) + \OO(h^2)\right) \\
			&= \Phi(y) + \OO(h),
		\end{aligned}
	\end{equation}
	which implies
	\begin{equation}\label{eq:DistanceProbDet}
		\begin{aligned}
			e^{h\diffL_h}\Phi(y) - \Phi(\Psi_h(y)) &= \frac{1}{2} Ch^{2p+1}\partial_{tt}\Psi_h(y) \cdot \nabla\Phi(y)\\
			&\quad +\frac{1}{2}Ch^{2p+1}\partial_t \Psi_h(y) \partial_t \Psi_h(y)^\top  \colon \nabla^2\Phi(y) + \OO(h^{2p+1}).
		\end{aligned}
	\end{equation}
	Let us remark that due to the smoothness of the flow we have
	\begin{equation}\label{eq:DistanceExactDet}
		e^{h\diffL}\Phi(y) - \Phi(\Psi_h(y)) = \OO(h^{q+1}).
	\end{equation}
	Combining \eqref{eq:DistanceExactDet} and \eqref{eq:DistanceProbDet} we have the one-step weak error of the probabilistic method on the original ODE, i.e., 
	\begin{equation}\label{eq:LocalWeakError}
		e^{h\diffL}\Phi(y) - e^{h\diffL_h}\Phi(y) = \OO(h^{\min\{2p+1, q+1\}}),
	\end{equation}
	which proves the desired result.
\end{proof}

\begin{remark} Let us remark that rigorously if $\partial_{tt}\Psi_h(y)$ is bounded independently of $y$ then the equality \eqref{eq:TaylorNumericalSolution} holds. In fact, as it can be noticed in \eqref{eq:TaylorNumericalSolution2}, a sufficient requirement is that $h^{p+1/2}\partial_{tt}\Psi_h(y)$ is bounded independently of $h$. 
\end{remark}

In order to obtain a result on the global order of convergence we need a further stability assumption, which is the same as Assumption 3 in \cite{CGS17}.
\begin{assumption}\label{as:Stability} The function $f$ and the distribution of the random time steps $H_k$, $k = 0, 1, \ldots$, are such that the operator $e^{h\diffL_h}$ satisfies for all functions $\psi\in \mathcal C^q_b(\R^d, \R)$ and a positive constant $L$, 
	\begin{equation}\label{eq:Stability}
		\sup_{u\in\R^d} \abs{e^{h\diffL_h}\psi(u)} \leq (1 + Lh)\sup_{u\in\R^d}\abs{\psi(u)},
	\end{equation}
	where $L$ may depend on $f$ and on the distribution of the random time steps, but not on $\psi$ or $h$.
\end{assumption}
\begin{remark} Let us remark that in order for $\Psi_h$ to satisfy \cref{as:PsiStrong}, i.e., for $\Psi_h$ to be of order $q$, the right hand side $f$ must be of class $\mathcal C_b^q(\R^d, \R^d)$ (see, e.g. \cite[Theorem II.3.1]{HNW93}). Therefore, in order to apply the bound \eqref{eq:Stability} to composite functions $\Phi \circ \phi_h \colon \R^d \to \R$ where $\Phi \in \mathcal C^\infty_b(\R^d, \R)$, by the chain rule we need \cref{as:Stability} to hold for functions in $C^q_b(\R^d, \R)$. This fact will be exploited in the proof of \cref{thm:weakOrder} below.
\end{remark}
We now give a lemma useful for bounding discrete sequences, which is taken from \cite[Lemma 1.6]{MiT04}.
\begin{lemma}\label{lem:RecurrenceBound} Suppose that for arbitrary $N$ and $k = 0, \ldots, N$ we have
		\begin{equation}
		e_k \leq (1 + Ah) e_{k-1} + Bh^r,
		\end{equation}
		where $h = T / N$, $A > 0$, $B \geq 0$, $r \geq 1$ and $e_k \geq 0$, $k = 0, \ldots, N$. Then
		\begin{equation}
		e_k \leq e^{AT}e_0 + \frac{B}{A}(e^{AT} - 1) h^{r-1}.
		\end{equation}
\end{lemma}

The proof of \cref{lem:RecurrenceBound} follows from the discrete Grönwall inequality. We can now state the main result on weak convergence. 

\begin{theorem}[Weak order]\label{thm:weakOrder} Let the assumptions of \cref{lem:WeakLocalOrder} and \cref{as:Stability} hold. Then, there exists a constant $C > 0$ independent of $h$ and of the initial condition such that for all functions $\Phi \in \mathcal C_b^l(\R^d,\R)$, with $l = \max\{q, 3\}$
	\begin{equation}\label{eq:weakOrder}
		\abs{\E\Phi(Y_k) - \Phi(y(kh)))} \leq Ch^{\min\{2p, q\}},
	\end{equation}
	for all $k = 1, 2, \ldots, N$ and $T = Nh$.
\end{theorem}
	
\begin{proof} Let us introduce the following notation
	\begin{equation}
		\begin{aligned}
			w_k(u) &= \Phi(\phi_{t_k}(u)),\\
			W_k(u) &= \E(\Phi(Y_k) \mid Y_0 = u).
		\end{aligned}
	\end{equation}
	By the triangle inequality and the Markov property \eqref{eq:MarkovProperty}, we have
	\begin{equation}
		\begin{aligned}
			\sup_{u\in\R^d}\abs{W_k(u) - w_k(u)} &\leq \sup_{u\in\R^d}\abs{e^{h\diffL}w_{k-1}(u) - e^{h\diffL_h}w_{k-1}(u)} \\
			&\quad + \sup_{u\in\R^d}\abs{e^{h\diffL_h}w_{k-1}(u) - e^{h\diffL_h}W_{k-1}(u)}.
		\end{aligned}
	\end{equation}
	We then apply \cref{lem:WeakLocalOrder} to the first term and \cref{as:Stability} to the second and denote $e_k \defeq \sup_{u\in\R^d}\abs{W_k(u) - w_k(u)}$, thus obtaining
	\begin{equation}
		e_k \leq Ch^{\min\{2p+1, q + 1\}} + (1 + Lh)e_{k-1}.
	\end{equation}
	We can therefore apply \cref{lem:RecurrenceBound} with $A = L$ and $r = \min\{2p+1, q + 1\}$, and therefore get for a constant $C > 0$
	\begin{equation}
		\sup_{u\in\R^d}\abs{w_k(u) - W_k(u)} \leq Ch^{\min\{2p, q\}},
	\end{equation}
	which is the desired result.
\end{proof}

\begin{remark}\label{rem:sde} In \cite{CGS17}, Conrad et al. define ordinary and stochastic modified equations in order to prove a result of weak convergence applying techniques of backward error analysis. In particular, they show that their probabilistic solver approximates in the weak sense a stochastic differential equation (SDE) where the deterministic part is given by the original ODE. For our probabilistic solver, it is possible to prove that the numerical solutions approximates in the weak sense the solution of an SDE which depends on the derivative of the map $t \mapsto \Psi_t(y)$. Such a construction is shown in the Appendix.
\end{remark}

\begin{remark} Let us recall that the random variable $Y_k$ given by RTS-RK is thought of as an approximation of $y(kh)$ regardless of the value of the sum of the random time steps. Hence, the comparison in \eqref{eq:weakOrder} is legitimate and does not induce time misalignment between true and numerical solutions. This basic property applies to all results in the following.
\end{remark}

\section{Mean square convergence analysis}\label{sec:StrongOrder}

The second property of \eqref{eq:RTS-RK} we analyze is its mean square order of convergence, which gives an indication on the path-wise distance between each realisation of the numerical solution and the exact solution of \eqref{eq:ODE}. Let us define the mean square order of convergence. 
\begin{definition} The numerical method \eqref{eq:RTS-RK} has mean square order of convergence $r$ for \eqref{eq:ODE} if there exists a constant $C > 0$ independent of $h$ and of the initial condition $y_0$ such that
	\begin{equation}
	\big(\E\norm{Y_k - y(kh)}^2\big)^{1/2} \leq Ch^r
	\end{equation}
	for all $k = 1, 2, \ldots, N$ and $T = Nh$.
\end{definition} 
\begin{remark} Let us remark that the mean square convergence is stronger than the traditional strong convergence, since, by Jensen's inequality 
	\begin{equation}
		\E{\norm{Y_k - y(kh)}} \leq \big(\E\norm{Y_k - y(kh)}^2\big)^{1/2} \leq Ch^r.
	\end{equation}	
\end{remark}

We start by analysing how the method converges with respect to the mean step size $h$ in the local sense, i.e., after one step of the numerical integration.
\begin{lemma}[Local mean square convergence]\label{thm:StrongOrderLocal} Under \cref{as:hStrong}, \cref{as:PsiStrong} and \cref{as:hImplicit} the numerical solution $Y_1$ given by one step of the RTS-RK method \eqref{eq:RTS-RK} satisfies 
	\begin{equation}\label{eq:StrongOrderLocal}
		\big(\E\norm{Y_1 - y(h)}^2\big)^{1/2} \leq C h^{\min\{p+1/2, q + 1\}},
	\end{equation}
	where $C$ is a real positive constant independent of $h$ and of the initial condition $y_0$ and the coefficients $p$, $q$ are given in the assumptions.
\end{lemma}
\begin{proof} By triangular and Young's inequalities we have for all $y \in \R^d$ 
	\begin{equation}
		\E\norm{\Psi_{H_0}(y) - \phi_h(y)}^2 \leq 2\E\norm{\Psi_{H_0}(y) - \Psi_h(y)}^2 + 2\norm{\Psi_h(y) - \phi_h(y)}^2.
	\end{equation}		
	We now consider \cref{as:PsiStrong} and \cref{as:hStrong}, thus getting
	\begin{equation}\label{eq:LocalErrorDecomposition}
	\begin{aligned}
		\E\norm{\Psi_{H_0}(y) - \phi_h(y)}^2 &\leq 2L_{\Psi}^2 \E\abs{H_0 - h}^2 + 2C_1 h^{2(q+1)} \\
		&= 2L_{\Psi}^2 C_2 h^{2p+1} + 2 C_1 h^{2(q+1)} \\
		&\leq C^2 h^{2\min\{p+1/2,q+1\}},
	\end{aligned}
	\end{equation}
	where $C_1$ and $C_2$ are the constants given in \cref{as:PsiStrong} and \cref{as:hStrong} respectively. This is the desired result with $C = \max\{2L_{\Psi}^2 C_2, 2 C_1\}^{1/2}$.
\end{proof}

As a consequence of the one-step convergence, we can prove a result of global mean square convergence.
\begin{theorem}[Global mean square convergence]\label{thm:StrongOrder} Let $f$ be globally Lipschitz and $t_k = kh$ for $k = 1, 2, \ldots, N$, where $Nh = T$. Then, under the assumptions of \cref{thm:StrongOrderLocal} the numerical solution given by \eqref{eq:RTS-RK} satisfies 
	\begin{equation}\label{eq:StrongGlobal}
		\sup_{k=1,2, \ldots, N} \big(\E\norm{Y_k - y(t_k)}^2\big)^{1/2} \leq C h^{\min\{p, q\}},
	\end{equation}
	where $C$ is a real positive constant independent of $h$ and of the initial condition.
\end{theorem}
In order to prove this result, let us introduce the following lemmas.
\begin{lemma}\label{lem:ODERepresentation} Given the ODE \eqref{eq:ODE} with $f$ globally Lipschitz, then for any $y$ and $w$ in $\R^d$ and $0 < h < 1$ we have
	\begin{align}
		&\norm{\phi_h(y) - \phi_h(w)} \leq (1 + Ch) \norm{y - w},\label{eq:FlowLipschitz}\\
		&\norm{\phi_h(y) - \phi_h(w) - (y - w)} \leq Ch\norm{y - w},\label{eq:FlowLipschitzBoundRem}
	\end{align}
	where $C$ is a positive constant independent of $h$ and of the initial condition $y_0$.
\end{lemma}

The proof of \cref{lem:ODERepresentation} follows from the global Lipschitz continuity of $f$ and the Grönwall inequality. We can now prove the main result on mean square convergence.
\begin{proof}[Proof of \cref{thm:StrongOrder}] In the following, we denote by $C$ a constant that does not depend on $h$ and on the initial condition $y_0$ whose value may change from line to line. Let us define $e_k^2 \defeq \E\norm{Y_k - y(t_k)}^2$. Adding and subtracting the exact flow applied to the numerical solution, we obtain
	\begin{equation}\label{eq:StrongErrorDecomp}
		\begin{aligned}
			e_{k+1}^2 = &\E\norm{\Psi_{H_k}(Y_k) - \phi_{h}(Y_k)}^2 + \E\norm{\phi_{h}(Y_k) - \phi_{h}(y(t_k))}^2 \\
					  &+ 2\E\Big(\big(\phi_{h}(Y_k) - \phi_{h}(y(t_k))\big)^\top \big(\Psi_{H_k}(Y_k) - \phi_{h}(Y_k)\big)\Big).
		\end{aligned}
	\end{equation}
	Let us consider the three terms in \eqref{eq:StrongErrorDecomp} separately. For the first term, we have by \cref{thm:StrongOrderLocal}
	\begin{equation}\label{eq:StrongFirstTerm}
		\E\norm{\Psi_{H_k}(Y_k) - \phi_{h}(Y_k)}^2 \leq Ch^{\min\{2p+1, 2(q+1)\}}.
	\end{equation}
    For the second term, due to \eqref{eq:FlowLipschitz}, we have
	\begin{equation}\label{eq:StrongSecondTerm}
		\E\norm{\phi_{h}(Y_k) - \phi_{h}(y(t_k))}^2 \leq (1 + Ch)^2 e_k^2.
	\end{equation}
	Let us now define $Z = \phi_{h}(Y_k) - \phi_{h}(y(t_k)) - (Y_k - y(t_k))$. Then we can rewrite the inner product as
	\begin{equation}\label{eq:StrongScalarProd}
	\begin{aligned}
		\E\Big(\big(\phi_{h}(Y_k) - \phi_{h}(y(t_k))\big)^\top \big(\Psi_{H_k}(Y_k) - \phi_{h}(Y_k)\big)\Big) = &\E\Big(\big(Y_k - y(t_k)\big)^\top \big(\Psi_{H_k}(Y_k) - \phi_{h}(Y_k)\big)\Big) \\
		&+ \E\Big(Z^\top \big(\Psi_{H_k}(Y_k) - \phi_{h}(Y_k)\big)\Big).
	\end{aligned}
	\end{equation}
	We bound the two terms in \eqref{eq:StrongScalarProd} separately. For the first term, by the law of total expectation, we have
	\begin{equation}
	\begin{aligned}
		\E\Big(\big(Y_k - y(t_k)\big)^\top \big(\Psi_{H_k}(Y_k) - \phi_{h}(Y_k)\big)\Big) &= \E\E\Big(\big(Y_k - y(t_k)\big)^\top \big(\Psi_{H_k}(Y_k) - \phi_{h}(Y_k)\big) \mid Y_k\Big)\\
		&= \E\Big(\big(Y_k - y(t_k)\big)^\top\E \big(\Psi_{H_k}(Y_k) - \phi_{h}(Y_k) \mid Y_k \big)\Big).
	\end{aligned}
	\end{equation}
	Applying Cauchy--Schwarz inequality to the outer expectation we get
	\begin{equation}
	\begin{aligned}
		\E\Big(\big(Y_k - y(t_k)\big)^\top \big(\Psi_{H_k}(Y_k) - \phi_{h}(Y_k)\big)\Big) &\leq \Big(\E\norm{\E\big(\Psi_{H_k}(Y_k) - \phi_{h}(Y_k) \mid Y_k\big)}^{2}\Big)^{1/2}e_k  \\
		&\leq  Ch^{\min\{2p+1, q+1\}}e_k,
	\end{aligned}
	\end{equation}
	where we applied \cref{lem:WeakLocalOrder}. We now consider the second term in \eqref{eq:StrongScalarProd}. By the Cauchy--Schwarz inequality we have
	\begin{equation}
		\E\Big(Z^\top \big(\Psi_{H_k}(Y_k) - \phi_{h}(Y_k)\big)\Big) \leq \big(\E\norm{Z}^{2}\big)^{1/2} \big(\E\norm{\Psi_{H_k}(Y_k) - \phi_{h}(Y_k)}^2\big)^{1/2}.
	\end{equation}
	We now apply \eqref{eq:FlowLipschitzBoundRem} and \cref{thm:StrongOrderLocal} to obtain
	\begin{equation}
		\E\Big(Z^\top \big(\Psi_{H_k}(Y_k) - \phi_{h}(Y_k)\big)\Big) \leq C h^{\min\{p+3/2, q+2\}} e_k.
	\end{equation}
	We can hence bound the scalar product in \eqref{eq:StrongScalarProd} with Young's inequality and assuming $h < 1$ as
	\begin{equation}\label{eq:StrongThirdTerm}
	\begin{aligned}
		\E\Big(\big(\phi_{h}(Y_k) - \phi_{h}(y(t_k))\big)^\top\big(\Psi_{H_k}(Y_k) - \phi_{h}(Y_k)\big)\Big) &\leq C h^{\min\{p+3/2, q+1\}} e_k \\
		&\leq \frac{he_k^2}{2} + C\frac{h^{\min\{2p+2,2q+1\}}}{2}.
	\end{aligned}
	\end{equation}
	Combining \eqref{eq:StrongFirstTerm}, \eqref{eq:StrongSecondTerm} and \eqref{eq:StrongThirdTerm}, we have
	\begin{equation}
		e_{k+1}^2 \leq Ch^{\min\{2p+1, 2q+1\}} + (1 + Ch)e_k^2,
	\end{equation}
	which implies the desired result by \cref{lem:RecurrenceBound} and since $e_0 = 0$.
\end{proof}

\begin{remark} Let us remark that the difference between global and local orders of convergence, i.e., between \eqref{eq:StrongOrderLocal} and \eqref{eq:StrongGlobal}, is not exactly one, as it usually is in the purely deterministic case. In fact, due to the independence of the random variables there is only a $1/2$ loss in the random part of the exponent, while the natural loss of one order is verified in the deterministic component.
\end{remark}
\begin{remark} As for the additive noise method proposed in \cite{CGS17}, the result of mean square convergence suggests that a reasonable choice for the noise scale $p$ is to fix $p=q$, where $q$ is the order of the Runge--Kutta method $\Psi_h$. In this way, the properties of convergence of the underlying deterministic method are preserved, while yielding a probabilistic interpretation of the numerical solution.
\end{remark}

\section{Mean square convergence of Monte Carlo estimators}\label{sec:MonteCarlo}

The third property we analyze is the mean-square convergence of Monte Carlo estimators drawn from the random time-stepping Runge--Kutta method. Let us consider a function $\Phi\in\mathcal{C}^\infty_b(\R^d, \R)$ with Lipschitz constant $L_\Phi$ and a final time $T > 0$. Moreover, let us introduce the notation $Z = \Phi(y(T))$ and $Z_N = \E\Phi(Y_N)$, where $N$ is such that $T = Nh$. In general, the quantity $Z_N$ is not accessible, and we have to replace it by its Monte Carlo estimator 
\begin{equation}\label{eq:MSE}
	\hat Z_{N, M} = M^{-1} \sksum_{i = 1}^M \Phi(Y_N^{(i)}).
\end{equation}
where $M$ is the number of realisations of the numerical solution and we denote by $\{Y_N^{(i)}\}_{i=1}^M$ a set of i.i.d. realisations of the numerical solution. Hence, we are interested in studying the mean square error of the Monte Carlo estimator, which is defined as
\begin{equation}
	\MSE(\hat Z_{N, M}) = \E(Z - \hat Z_{N, M})^2.
\end{equation}
In the following result, we prove that this quantity converges to zero independently of the number of trajectories $M$, in the limit $h \to 0$.

\begin{theorem}\label{thm:MSEMonteCarlo} Under \cref{as:hStrong}, \cref{as:hImplicit} and \cref{as:PsiStrong}, the Monte Carlo estimator $\hat Z_{N,M}$ satisfies
	\begin{equation}\label{eq:MSEBound}
		\MSE(\hat Z_{N,M}) \leq C\Big(h^{2\min\{2p, q\}} + \frac{h^{2\min\{p, q\}}}{M}\Big),
	\end{equation}
	where $C$ is a positive constant independent of $h$ and $M$.
\end{theorem}
\begin{proof} Thanks to the classic decomposition of the MSE, we have
	\begin{equation}\label{eq:MSEDecomposition}
		\MSE(\hat Z_{N,M}) = \Var \hat Z_{N,M}  + \big(\E(\hat Z_{N,M} - Z)\big)^2.
	\end{equation}
	Due to the unbiasedness of the Monte Carlo estimator $\hat Z_{N, M}$ and applying \cref{thm:weakOrder} to the second term, we have
	\begin{equation}\label{eq:MSEWeakOrder}
		\MSE(\hat Z_{N,M}) \leq \Var \hat Z_{N,M}  + Ch^{2\min\{2p, q\}}.
	\end{equation}
	The variance of the estimator can be trivially bounded by exploiting the Lipschitz continuity of $\Phi$ and the independence of the samples {as}
	\begin{equation}\label{eq:MSELipschitz}
	\begin{aligned}
		\Var\hat Z_{N, M} &= M^{-1} \Var\big(\Phi(Y_N)\big)\\
		&\leq M^{-1} \E\big(\Phi(Y_N) - \Phi(y(T))\big)^2 \\
		&\leq M^{-1} L_\Phi^2 \E \norm{Y_N - y(T)}^2.
	\end{aligned}
	\end{equation}
	Applying \cref{thm:StrongOrder} we get
	\begin{equation}
		\Var\hat Z_{N,M} \leq M^{-1} L_\Phi^2 Ch^{2\min\{p, q\}},
	\end{equation}
	which proves the desired result.
\end{proof}
Let us remark that with the choice $p = q$, which is the minimum $p$ for which the order of convergence of the underlying deterministic method is not affected by the probabilistic setting, we have $\MSE(\hat Z_{N,M}) \leq Ch^{2q}$ with $M = 1$. Hence, the Monte Carlo estimators drawn from \eqref{eq:RTS-RK} converge in the mean square sense independently of the number of samples $M$ in \eqref{eq:MSE}. In the sub-optimal case $p < q$, one should carefully select the number of trajectories $M$ so that the two terms in \eqref{eq:MSEBound} are balanced. In particular, this would lead to
\begin{equation}
	M = \begin{cases} \OO(1), & \mbox{if } p \geq q, \\
	\OO(h^{2(p-q)}), & \mbox{if } p < q \leq 2p,\\
	\OO(h^{-2p}), & \mbox{if } 2p < q,
	\end{cases}
\end{equation}
where the notation $M=\OO(h^r)$ for a real number $r$ means that there exist constants $C_1$ and $C_2$ such that $C_1 h^r \leq M \leq C_2 h^r$.
\begin{remark} Let us remark that in order to have uncertainty quantification for a fixed value $h > 0$ it is necessary to draw a sample with $M > 1$, since otherwise the probability distribution over the numerical solution would be a {Dirac} delta. \cref{thm:MSEMonteCarlo} does not provide an indication of how the value of $M$ should be chosen in order to have a good empirical description of the probability measure induced by the RTS-RK method, but still ensures quantitatively that the Monte Carlo estimators drawn from this distribution have a good quality. 
\end{remark}

\section{Conservation of first integrals}\label{sec:GeomProperties}
Numerical methods for ODEs are often studied in terms of their geometric properties \cite{HLW06}. In particular, we investigate here whether the random choice of time steps in \eqref{eq:RTS-RK} spoils the properties of the underlying deterministic Runge--Kutta method. Let us recall the definition of first integral for an ODE.
\begin{definition} Given a function $I\colon\R^d\to\R$, then $I(y)$ is a first integral of \eqref{eq:ODE} if $I'(y)f(y) = 0$ for all $y \in \R^d$. 
\end{definition}	
If this property of the ODE is conserved by a numerical integrator, i.e., if for the any $y\in\R^d$ it is true that $I(\Psi_h(y)) = I(y)$, then we say that the numerical method conserves the first integral. In particular, this implies that the first integral $I$ is conserved along the trajectory of the numerical solution, i.e., $I(y_k) = I(y_0)$ for all $k\geq 0$.
	
\begin{example} To illustrate this concept we first discuss the case of linear first integrals, which can be seen as a general case of the conservation of mass in physical systems. Let us consider a linear first integral $I(y) = v^\top y$ and any Runge--Kutta method with coefficients $\{b_i\}_{i=1}^s$, $\{a_{ij}\}_{i,j=1}^s$. Then, we have for a time step $H_0 > 0$
	\begin{equation}
		I(Y_1) = v^\top y_0 + H_0 \sksum_{i=1}^s b_iv^\top f(y_0 + H_0\sksum_{j=1}^{s} a_{ij}K_j),
	\end{equation}
	where $\{K_i\}_{i=1}^s$ are the internal stages of the Runge--Kutta method. Since $I(y)$ is a first integral, $v^\top f(y) = 0$ for any $y \in \R^d$. Hence $I(Y_1)  = I(y_0)$ and iteratively $I(Y_k) = I(y_0)$ for all $k \geq 0$ along the numerical trajectory. The equality above shows that any RTS-RK method conserves linear first integrals path-wise, or in the strong sense. 
\end{example}

It is known that no Runge--Kutta method can conserve any polynomial invariant of order $n \geq 3$ \cite[Theorem IV.3.3]{HLW06}. Nonetheless, for some particular problems there exist tailored Runge--Kutta methods which can conserve polynomial invariants of higher order. We therefore can state the following general result.
\begin{theorem}\label{thm:PolyInvariants} Let $I(y)$ be a first integral for \eqref{eq:ODE} and $\Psi_h$ be the numerical flow of a Runge--Kutta scheme for \eqref{eq:ODE}. If the scheme defined by $\Psi_h$ conserves $I(y)$ for any $h > 0$, then the numerical method \eqref{eq:RTS-RK} conserves $I(y)$ almost surely.
\end{theorem}
\begin{proof} If $I(\Psi_h(y)) = I(y)$ for any $h$, then $I(\Psi_{H_0}(y)) = I(y)$ almost surely for any value that $H_0$ can assume.
\end{proof}

We now consider quadratic first integrals, i.e., first integrals of the form $I(y) = y^\top S y$ with $S$ a symmetric matrix, which are conserved by Runge--Kutta methods that satisfy the hypotheses of Cooper's theorem \cite[Theorem IV.2.2]{HLW06}. The conservation of quadratic first invariants is of the utmost importance, e.g., for Hamiltonian systems, as it implies the symplecticity of the scheme. It is known \cite[Theorem IV.2.1]{HLW06} that all Gauss methods conserve quadratic first integrals. The simplest member of this class of methods is the implicit midpoint rule, which is a one-stage method defined by coefficients $b_1 = 1$ and $a_{11} = 1/2$.
\begin{corollary}\label{thm:QuadraticInvariants} If the Runge--Kutta scheme defined by $\Psi_h$ conserves quadratic first integrals then the numerical method \eqref{eq:RTS-RK} conserves quadratic first integrals almost surely.
\end{corollary}
\begin{proof} This result is a direct consequence of \cref{thm:PolyInvariants}. \end{proof}

The properties above for the RTS-RK method are not satisfied by the additive noise method presented in \cite{CGS17}. In particular, let us remark that the conservation of first integrals is exact for any trajectory of the RTS-RK method, and is not an average property. In other words, we can say that \eqref{eq:RTS-RK} conserves linear first integrals in the strong sense. For the additive noise numerical method \eqref{eq:AN-RK}, we have
\begin{equation}
	\begin{aligned}
	I(Y_1) &= v^\top  y_0 + h \sksum_{i=1}^s b_iv^\top  f(y_0 + h\sksum_{j=1}^{s} a_{ij}K_j) + v^\top  \xi_0(h), \\
	&= v^\top  (y_0 + \xi_0(h)).
	\end{aligned}
\end{equation}
If the random variable $\xi_0$ is zero-mean, then $\E I(Y_1) = I(y_0)$ and iteratively along the solution $\E I(Y_k) = I(y_0)$. Linear first integrals are therefore conserved in average, but not in a path-wise fashion.

For quadratic first integrals, we have instead that the additive noise method does not conserve them neither path-wise nor in the weak sense, as we have
\begin{equation}
\begin{aligned}
	I(Y_1) &= (\Psi_h(y_0) + \xi_0(h))^\top  S (\Psi_h(y_0) + \xi_0(h)) \\
	&= I(y_0) + 2\xi_0(h)^\top  S  \Psi_h(y_0) + \xi_0(h)^\top  S \xi_0(h).
\end{aligned}
\end{equation}
If the random variables are zero-mean and if there exists a matrix $Q$ such that $\E\xi_0(h)\xi_0(h)^\top  = Qh^{2p + 1}$ for some $p \geq 1$ (Assumption 1 in \cite{CGS17}) we then have
\begin{equation}\label{eq:BiasQuadraticAddNoise}
\E I(Y_1) = I(y_0) + Q : S h^{2p + 1}.
\end{equation}
Hence, along the trajectories of the solution a bias is introduced in the first integral which persists even in the mean sense. In general, \cref{thm:PolyInvariants} is not valid for the additive noise method, as the random contribution drives the first integral far from its true value at each time step. In practice, this could produce large deviations of the numerical approximation from the true solution, especially in the long time regime.

\section{Hamiltonian systems}\label{sec:Hamiltonian} A class of dynamical systems of particular interest for their geometric properties is the class of Hamiltonian systems. Given a function $Q \colon \R^{2d} \to \R$, called the Hamiltonian, Hamiltonian systems can be written as
\begin{equation}\label{eq:ODEHam}
y' = J^{-1}\nabla Q(y), \quad y(0) = y_0 \in \R^{2d},
\end{equation}
where the matrix $J\in\R^{2d \times 2d}$ is defined as
\begin{equation}
J = \begin{pmatrix} 0 & I \\ -I & 0 \end{pmatrix},
\end{equation}
and where $I$ is the identity matrix in $\R^{d\times d}$. The Hamiltonian $Q$ is a first integral for \eqref{eq:ODEHam}, hence we require numerical integrators to conserve the energy, or at least not to deviate from its true value in an uncontrolled fashion. As it was shown in the previous section, when $Q$ is a polynomial it is possible to obtain exact conservation with deterministic integrators and with their probabilistic counterparts obtained with the RTS-RK method. If $Q$ is not a polynomial, exact conservation is in general not achievable, but a good approximation of the energy over long time spans is achievable through the notion of symplectic differentiable maps.

\begin{definition}[Definition VI.2.2 in \cite{HLW06}] Let $U\subset \R^{2d}$ be a non-empty open set. A differentiable map $g\colon U \to \R^{2d}$ is called symplectic if the Jacobian matrix $g'$ is everywhere symplectic, i.e., if
	\begin{equation}\label{eq:SymplecticMap}
	(g')^\top J g' = J.
	\end{equation}	
\end{definition}
It is well-known that the flow $\phi_t\colon\R^{2d}\to\R^{2d}$ of any system of the form \eqref{eq:ODEHam} is symplectic. In a natural manner, a numerical integrator is called symplectic if its numerical flow $\Psi_h$ is a symplectic map {whenever it is applied to a smooth Hamiltonian system \cite[Definition VI.3.1]{HLW06}}. In the following, we will analyse both the local and global properties of the RTS-RK method built on symplectic integrators and applied to \eqref{eq:ODEHam}.

\subsection{Symplecticity of the RTS-RK method} It has been pointed out \cite[Section VIII.1]{HLW06} that applying an adaptive step size technique to a symplectic method can destroy its symplecticity. Therefore, Skeel and Gear \cite{SkG92} write any adaptive technique in terms of a map $\tau(y, h)$ such that the $k$-th time step $h_k$ is selected as $h_k = \tau(y_k, h)$, where $h$ is a base value for the time step. Hence, in order to have again a symplectic method {for variable time steps}, the new condition to be satisfied is
\begin{equation}\label{eq:SymplecticCondition}
V^\top J V = J, \quad V = \pdv{\Psi_{\tau(y, h)}(y)}{y} + \pdv{\Psi_{\tau(y, h)}(y)}{t}\pdv{\tau(y, h)}{y}^\top.
\end{equation}
Let us now consider the RTS-RK method based on a symplectic deterministic integrator. We have the following lemma. 
\begin{lemma}\label{lem:SympRTSRK} If the flow $\Psi_h$ of the deterministic integrator is symplectic, then the flow of the random time-stepping probabilistic method \eqref{eq:RTS-RK} is symplectic.
\end{lemma}
\begin{proof} For the RTS-RK scheme, the $k$-th time step $H_k$ is generated by a random mapping as $H_k = \tau(y, h) = \tau(h) = h\Theta_k$, where $\Theta_k$ are appropriately scaled random variables such that $H_k$ satisfies \cref{as:hStrong}. Hence, $\tau$ is independent of $y$, i.e., $\partial_y \tau(y, h) = 0$, and with the notation introduced above 
	\begin{equation}
	V = \pdv{\Psi_{\tau(h)}(y)}{y}.
	\end{equation}
	Therefore, by the symplecticity of $\Psi_t$ the condition $V^\top J V = J$ is satisfied and the flow map of the RTS-RK method is symplectic.
\end{proof}
Let us remark that the local symplecticity of the flow map is not sufficient for good conservation of the Hamiltonian for the numerical solution. Global properties of approximation of the energy are therefore presented below.

\subsection{Long-time conservation of Hamiltonians}\label{sec:Hamiltonian_2} We now wish to study the mean conservation of the Hamiltonian along the trajectories of the RTS-RK method based on symplectic integrators. Our goal is obtaining a bound on the quantity $\E\abs{Q(Y_n) - Q(y_0)}$ that holds over long times. Showing theoretically long time conservation of the energy function in Hamiltonian systems requires backward error analysis. In the following, we will introduce the basis of this technique and show how they apply to our probabilistic integrator. {For further details, a comprehensive treatment of backward error analysis ought to be found in \cite[Chapter IX]{HLW06}.}

The first ingredient needed to perform a rigorous backward error analysis is a rather strong assumption on the regularity of the ODE, see e.g. \cite[Section IX.7]{HLW06}.
\begin{assumption}\label{as:RegHamiltonian} The function $f$ is analytic in a neighbourhood of the initial condition $y_0$ and there exist constants $C, R > 0$ such that $\norm{f(y)} \leq C$ for $\norm{y - y_0} \leq 2R$.
\end{assumption}
In general, backward error analysis is based on determining a modified equation $y' = \tilde f(y)$ such that the numerical approximation is its exact solution. Hence, the function $\tilde f$ will both depend on the original ODE and on the numerical flow map $\Psi_h$. In particular, for an integrator of order $q$ the modified equation is given by a function $\tilde f$ defined as
\begin{equation}
\tilde f(y) = f(y) + h^q f_{q+1}(y) + h^{q+1} f_{q+2}(y) + \ldots,
\end{equation}
where the functions $\{f_i\}_{i>q}$ are uniquely determined by $f$, its derivatives and by the coefficients of the Runge--Kutta method. The exactness of the numerical solution for the modified equation is nonetheless only formal, as the infinite sum defining $\tilde f$ is not guaranteed to converge. Thus, it is necessary to truncate the sum in order to perform a rigorous analysis, i.e.,
\begin{equation}\label{eq:ModifiedRHS}
\tilde f(y) = f(y) + h^q f_{q+1}(y) + h^{q+1} f_{q+2}(y) + \ldots + h^{N-1}f_N(y).
\end{equation}
where $q < N < \infty$ is the truncation index. Let us remark that in the following we will always refer to the truncated function above when using the symbol $\tilde f$. The truncation of the infinite sum implies that the numerical solution is not exact for the modified equation anymore. In particular, the error committed over one step on the modified equation is given by (see e.g. \cite[Theorem IX.7.6]{HLW06})
\begin{equation}\label{eq:OneStepModEq}
\norm{\tilde \phi_h(y) - \Psi_h(y)} \leq Che^{-\kappa/h},
\end{equation}
where $\tilde \phi$ is the exact flow of the modified equation and $\kappa$ and $C$ are constants depending on the coefficients of the method and on the regularity of $f$.

It is possible to prove (see e.g. \cite[Section IX.8]{HLW06}) that for a Hamiltonian system \eqref{eq:ODEHam} and a symplectic integrator the modified equation is still a Hamiltonian system, i.e., there exists a modified Hamiltonian $\tilde Q$ defined as
\begin{equation}\label{eq:ModifiedHamiltonianTrunc}
\tilde Q(y) = Q(y) + h^q Q_{q+1}(y) + \ldots + h^{N-1} Q_N(y),
\end{equation}
such that $\tilde f = J^{-1} \nabla \tilde Q$. The estimate \eqref{eq:OneStepModEq} implies that the modified Hamiltonian is almost conserved by the symplectic integrator. In particular, if $Q$ is Lipschitz, we have
\begin{equation}\label{eq:LocalErrorHamiltonianDet}
	\abs{\tilde Q(\Psi_h(y)) - \tilde Q(y)} \leq Che^{-\kappa/h}.
\end{equation}
The bound above guarantees that the modified Hamiltonian is well approximated for a long time, and as a consequence that the original Hamiltonian is almost conserved for the same time span. In particular, the following result is valid, see e.g. \cite[Theorem IX.8.1.]{HLW06} or \cite{BeG94}.
\begin{theorem}\label{thm:DetHamiltonian} Under \cref{as:RegHamiltonian} and for $h$ sufficiently small, if the numerical solution $y_n$ given by a symplectic method of order $q$ applied to an Hamiltonian system is close enough to the initial condition $y_0$, then 
	\begin{equation}
	\begin{aligned}
	\tilde Q(y_n) &= \tilde Q(y_0) + \OO(e^{-\kappa / 2h}), \\
	Q(y_n) &= Q(y_0) + \OO(h^q).
	\end{aligned}
	\end{equation}
	over exponentially long time intervals $nh \leq e^{\kappa / 2h}$.
\end{theorem}
The randomisation of the time steps implies that a general modified equation does not exist. Nonetheless, due to \cref{lem:SympRTSRK}, it is possible to construct locally a random Hamiltonian modified equation at each time step. We thus define at each step the random modified Hamiltonian as
\begin{equation}\label{eq:ModifiedHamiltonianTruncRand}
\hat Q_j(y) = Q(y) + H_j^q Q_{q+1}(y) + \ldots + H_j^{N-1} Q_N(y).
\end{equation}
As for the deterministic case, the random modified Hamiltonian $\hat Q$ will be almost conserved by the numerical flow. In particular, we define the random local truncation error as
\begin{equation}\label{eq:etaJ}
	\eta_j \defeq \hat Q_j(\Psi_{H_j}(y)) - \hat Q_j(y),
\end{equation} 
which, in light of \eqref{eq:LocalErrorHamiltonianDet}, satisfy 
\begin{equation}\label{eq:BoundEtaJ}
\abs{\eta_j} \leq C H_j e^{-\kappa / H_j},
\end{equation}
almost surely. In order to prove the conservation of the Hamiltonian over long time for the RTS-RK method, it is necessary to introduce a technical assumption on the higher moments of the random time steps.
\begin{assumption}\label{as:RTSMoments} There exists $\bar r > 1$ such that for any $1 < r < \bar r$, the random time steps $\{H_j\}_{j\geq 0}$ satisfy
	\begin{equation}
	\E H_j^r = h^r + C_rh^{2p+r-1},
	\end{equation}
	where $p$ is defined in \cref{as:hStrong} and $C_r > 0$ satisfies $C_{2r} > 2C_r$ and is independent of $h$. Moreover, there exists $m, M > 0$ with $M > m$ such that $mh \leq H_j \leq Mh$ almost surely for all $j \geq 0$.
\end{assumption}
This assumption guarantees that the higher moments of the random time steps are close to the corresponding powers of $h$ in the mean and mean square sense. In particular, it is possible to verify that
\begin{equation}
\begin{aligned}
	&\E(H_j^r - h^r) = C_r h^{2p+r-1},\\
	&\E(H_j^r - h^r)^2 = (C_{2r} - 2C_r) h^{2p+2r-1}.
\end{aligned}
\end{equation}
Then, for any $r, s > 1$ such that $r + s < R$, it holds
\begin{equation}\label{eq:higherHigherMoments}
\begin{aligned}
	&\E(H_j^{r+s} - h^{r+s}) = \hat C_{r,s} h^s \E(H_j^r - h^r),\\
	&\E(H_j^{r+s} - h^{r+s})^2 = \tilde C_{r,s} h^{2s} \E(H_j^r - h^r)^2,
\end{aligned}
\end{equation}
where $\hat C_{r,s} = C_{r+s}/C_r$ and $\tilde C_{r,s} = (C_{2(r+s)} - 2C_{r+s})/(C_{2r} - 2C_r)$. Finally, let us remark that \cref{as:RTSMoments} is satisfied for the uniform random time steps $H_j \iid \mathcal U(h-h^{p+1/2}, h+h^{p+1/2})$ introduced in \cref{ex:uniformH}. Let us now prove a bound on the random variables $\eta_j$ defined in \eqref{eq:etaJ}.
\begin{lemma}\label{lem:BoundEtaJ} {Suppose that \cref{as:hStrong}, \cref{as:hImplicit} and \cref{as:RTSMoments} hold true, and suppose that $0 < h \leq 1$. Then} the random variables $\eta_j$ satisfy
	\begin{equation}
	\E \abs{\eta_j}^r \leq C h^{\min\{r, p+r-3/2\}} e^{-r\kappa/(Mh)}, 
	\end{equation}
	where $C > 0$ is independent of $h$ and for all $r \in \N$ with $r \geq 1$.
\end{lemma}
\begin{proof} The proof is given in the Appendix. \end{proof}
Let us furthermore introduce two lemmas, which will be employed for proving long-time conservation of Hamiltonians. {Let us remark that in \cref{lem:Expansion} the values $n, q, N$ indicate generic positive integers.}   
\begin{lemma}\label{lem:Expansion}  Let $n, q, N$ be positive integers with $N > q$, and let us define the sets of real numbers $a = a_{n,q,N} \defeq \{a_{jk}, j = 0, \ldots, n-1$, $k = q, \ldots, N-1\}$ and $b = b_n \defeq \{b_j, j = 0, \ldots, n-1\}$. Then
	\begin{equation}
	\Big(\sum_{j=0}^{n-1} \Big(\sum_{k=q}^{N-1} a_{jk} + b_j\Big)\Big)^2 = \sum_{j=0}^{n-1} a_{jq}^2 + 2 \sum_{j=1}^{n-1} \sum_{i=0}^{j-1} a_{jq}a_{iq} + R(a) + S(a, b),
	\end{equation}
	where the remainder $R(a)$ can be written as $R = R_1 + R_2 + R_3$, with
	\begin{equation}
	\begin{alignedat}{2}
		&R_1(a) = \sum_{j=0}^{n-1}\sum_{k=q+1}^{N-1} a_{jk}^2,	&&R_2(a) = 2 \sum_{j=0}^{n-1}\sum_{k=q+1}^{N-1}\sum_{l=q}^{k-1} a_{jk} a_{jl},\\
		&R_3(a) = 2 \sum_{j=1}^{n-1} \sum_{i=0}^{j-1} \sum_{k=q}^{N-1}\sum_{\substack{l=q \\ l+k > 2q}}^{N-1} a_{jk}a_{il},
	\end{alignedat}
	\end{equation}
	and the remainder $S(a, b)$ can be written as $S = S_1 + S_2 + S_3 + S_4$, with
	\begin{equation}
	\begin{alignedat}{2}
		&S_1(a, b) = \sum_{j=0}^{n-1} b_j^2,  &&S_2(a, b) = 2\sum_{j=1}^{n-1} \sum_{i=0}^{j-1} b_i b_j, \\
		&S_3(a, b) = 2\sum_{j=1}^{n-1} \sum_{k=q}^{N-1} b_j a_{jk}, \qquad &&S_4(a, b) = 2\sum_{j=1}^{n-1}\sum_{i=0}^{n-1}\Big(b_j \sum_{k=q}^{N-1} a_{ik} + b_i \sum_{k=q}^{N-1} a_{jk}\Big).\\
	\end{alignedat}
	\end{equation}
\end{lemma}
\begin{proof} The proof is given in the Appendix. \end{proof}

\begin{lemma}\label{lem:Remainder} Let \cref{as:hStrong} hold with $p \geq 3/2$ and $h < 1$, and let \cref{as:PsiStrong}, \cref{as:RegHamiltonian} and \cref{as:RTSMoments} hold. {Moreover, let $q$ be specified in \cref{as:PsiStrong} and $N$ be the truncation index of the modified right hand side \eqref{eq:ModifiedRHS}}. Let us consider the sets of real-valued random variables $\Delta \defeq \{\Delta_{j,k}(H_j^k - h^k), j = 0, \ldots, n-1$, $k = q, \ldots, N-1\}$, where $\Delta_{j,k} \defeq  Q_{k+1}(Y_j) - Q_{k+1}(Y_{j+1})$ and $\eta \defeq \{\eta_j, j = 0, \ldots, n-1\}$.  Then, with the notation of \cref{lem:Expansion}, there exist positive constants $C_1$, $C_2$ independent of $h$ and $n$, but possibly dependent on $q$ and $N$, such that
	\begin{equation}
	\begin{aligned}
		\E R(\Delta) &\leq C_1 \big(t_n h^{2(p + q + 1/2)} + t_n^2 h^{2(2p + q - 1/2)}\Big),\\
		\E S(\Delta, \eta) &\leq C_2 \big( (t_n h + t_n^2) e^{-2\kappa/(Mh)} + (t_n h^{p+q+1/2} + t_n^2 h^{2p+q-1})e^{-\kappa/(Mh)}\big),
	\end{aligned}
	\end{equation}
	where $t_n = nh$.
\end{lemma}
\begin{proof} The proof is given in the Appendix. \end{proof}
It is now possible to prove a result of long conservation of the Hamiltonian for symplectic RTS-RK methods.
\begin{theorem}\label{thm:RTSHamiltonian} Let $0 < h \leq 1$. {Suppose that \cref{as:hStrong} holds for $p \geq 3/2$, that \cref{as:hImplicit} and \cref{as:RegHamiltonian} hold, and that \cref{as:RTSMoments} holds with $\bar r$ sufficiently large. Moreover, let $Y_n$ be the solution given by the RTS-RK method built on a symplectic integrator of order $q$ applied to a Hamiltonian system with Hamiltonian $Q$. If $Y_0 = y_0$ and the numerical solution $Y_n$ is close enough to the initial condition $y_0$ almost surely, then there exist a constant $C > 0$ independent of $h$ and $n$ such that}
		\begin{equation}
		\E \abs{Q(Y_n) - Q(y_0)} \leq C h^q,
		\end{equation}
		for time intervals of length 
		\begin{equation}
		t_n = \OO\big(\min\{h^{1-2p}, e^{\kappa/(4Mh)} h^{-(2p+2q-1)/4}, e^{\kappa/(2Mh)}\}\big)
		\end{equation}
		where $p$ is given in \cref{as:hStrong} and $M$ in \cref{as:RTSMoments}.
\end{theorem}
\begin{proof} In the following proof, we denote by $C$ a positive constant independent of $h$ and $n$ which can possibly change value from line to line. Let us first consider the modified Hamiltonian $\tilde Q$ and expand the difference $\tilde Q(Y_n) - \tilde Q(y_0)$ in a telescopic sum as
	\begin{equation}\label{eq:ProofHamiltionianBase}
	\tilde Q(Y_n) - \tilde Q(y_0) = \sum_{j=0}^{n-1} \big(\tilde Q(Y_{j+1}) - \tilde Q(Y_j)\big).
	\end{equation}
	We then consider each element of the sum, add and subtract the random modified Hamiltonian $\hat Q_j$ computed in $Y_{j+1}$ thus obtaining
	\begin{equation}
	\begin{aligned}
	\tilde Q(Y_{j+1}) - \tilde Q(Y_j) &= \tilde Q(Y_{j+1}) - \hat Q_j(Y_{j+1}) + \hat Q_j(Y_{j+1}) - \tilde Q(Y_j) \\
	&= \tilde Q(Y_{j+1}) - \hat Q_j(Y_{j+1}) + \hat Q_j(Y_j) - \tilde Q(Y_j) + \eta_j.
	\end{aligned}
	\end{equation}
	Hence, {by applying the definition \eqref{eq:ModifiedHamiltonianTrunc} of $\tilde Q$ and \eqref{eq:ModifiedHamiltonianTruncRand} of $\hat Q_j$}, we get
	\begin{equation}
	\tilde Q(Y_{j+1}) - \tilde Q(Y_j) = \sum_{k=q}^{N-1} (H_j^k - h^k)\Delta_{j,k} + \eta_j,
	\end{equation}
	where $\Delta_{j,k}$ is defined in \cref{lem:Remainder}.	Going back to \eqref{eq:ProofHamiltionianBase}, applying Jensen's inequality and \cref{lem:Expansion} we obtain
	\begin{equation}\label{eq:ProofHamiltonianDecomp}
	\begin{aligned}
	\big(\E \abs{\tilde Q(Y_n) - \tilde Q(y_0)}\big)^2 &\leq \E \Big(\sum_{j=0}^{n-1} \Big(\sum_{k=q}^{N-1}(H_j^k - h^k)\Delta_{j,k} + \eta_j \Big) \Big)^2 \\
	&= \sum_{j=0}^{n-1}\E\big((H_j^q - h^q)^2 \Delta_{j,q}^2\big) \\
	&\quad + 2\sum_{j=1}^{n-1}\sum_{i=0}^{j-1} \E\big((H_j^q - h^q)\Delta_{j,q}(H_i^q - h^q) \Delta_{i,q}\big) + \E R(\Delta) + \E S(\Delta, \eta).
	\end{aligned}
	\end{equation}
	The first term in \eqref{eq:ProofHamiltonianDecomp} satisfies
	\begin{equation}\label{eq:BoundFirstTerm}
		\Big(\sum_{j=0}^{n-1}\E\big((H_j^q - h^q)^2 \Delta_{j,q}^2\big)\Big)^{1/2} \leq C \sqrt{t_n} h^{p+q},
	\end{equation}
	due to \eqref{eq:RemainderFirstBound}. Now, considering {\eqref{eq:RemainderFourthBound}}, we obtain that the second term in \eqref{eq:ProofHamiltonianDecomp} satisfies
	\begin{equation}
		\Big(2\sum_{j=1}^{n-1}\sum_{i=0}^{j-1} \E\big((H_j^q - h^q)\Delta_{j,q}(H_i^q - h^q) \Delta_{i,q}\big)\Big)^{1/2} \leq  C t_n h^{2p+q-1}.
	\end{equation}
	For the remainder term $\E R(\Delta)$, due to \cref{lem:Remainder} we get
	\begin{equation}
		(\E R(\Delta))^{1/2} \leq C \big(\sqrt{t_n} h^{p + q + 1/2} + t_n h^{2p + q - 1/2}\Big).
	\end{equation}
	For the remainder term $\E S(\Delta, \eta)$, due to \cref{lem:Remainder} and since $h \leq 1$ and $p \geq 3/2$ by assumption, we get
	\begin{equation}\label{eq:BoundS}
	\begin{aligned}
		(\E S(\Delta, \eta))^{1/2} &\leq C \Big( t_n^2 \big(e^{-2\kappa/(Mh)} + h^{p+q+1/2}e^{-\kappa/(Mh)}\big) \Big)^{1/2}\\
		&\leq C t_n \big(e^{-\kappa/(Mh)} + h^{(2p+2q+1)/4}e^{-\kappa/(2Mh)}\big).
	\end{aligned}
	\end{equation}
	Finally, taking the square root of both sides of \eqref{eq:ProofHamiltonianDecomp}, replacing the expressions we obtained above and since $h \leq 1$, we get that the modified Hamiltonian satisfies
	\begin{equation}
	\begin{aligned}
		\E \abs{\tilde Q(Y_n) - \tilde Q(y_0)} \leq C \Big(\sqrt{t_n} h^{p+q} + t_n h^{2p+q-1} + t_n \big(e^{-\kappa/(Mh)} + h^{(2p+2q+1)/4}e^{-\kappa/(2Mh)}\big)\Big).
	\end{aligned}
	\end{equation}
	Hence, imposing for a constant $C > 0$
	\begin{equation}
		t_n \leq C \min\{h^{1-2p}, e^{\kappa/(4Mh)} h^{-(2p+2q-1)/4}, e^{\kappa/(2Mh)}\},
	\end{equation}
	and since exponential terms are dominated by polynomial terms (see e.g. \cite[Theorem IX.8.1]{HLW06}), we obtain
	\begin{equation}\label{eq:HamiltonianTheoremSemiResult}
		\E \abs{\tilde Q(Y_n) - \tilde Q(y_0)} \leq Ch^q.
	\end{equation}
	Finally, applying the triangle inequality, since {for all $y \in \R^d$ it holds} $\abs{Q(y) - \tilde Q(y)} \leq Ch^q$ by definition of the modified Hamiltonian $\tilde Q$ {and due to \eqref{eq:HamiltonianTheoremSemiResult}} we get
	\begin{equation}
	\begin{aligned}
		\E\abs{Q(Y_n) - Q(y_0)} &\leq \E\abs{Q(Y_n) - \tilde Q(Y_n)} + \E \abs{Q(y_0) - \tilde Q (y_0)} + \E\abs{\tilde Q(Y_n) - \tilde Q(y_0)} \\
		&\leq C h^q,
	\end{aligned}
	\end{equation}
	which is the desired result.
\end{proof}

\begin{remark} The result of \cref{thm:RTSHamiltonian} is consistent with the theory of deterministic symplectic integrators. In fact, in the limit $p \to \infty$, one can choose the coefficient $M$ in \cref{as:RTSMoments} arbitrarily close to $1$ and we have 
	\begin{equation}
	\E \abs{Q(Y_n) - Q(y_0)} = \OO(h^q),
	\end{equation}
	for exponentially long time spans $t_n = \OO\big(e^{\kappa/(2h)}\big)$, which is consistent with the theory of deterministic symplectic integrators summarised by \cref{thm:DetHamiltonian}.
\end{remark}
\begin{remark} It has been observed (see for example \cite{Hai97, HLW06}) that adopting variable step sizes in symplectic integration destroys the good properties of conservation of the Hamiltonian. In particular, the error on the Hamiltonian has a linear drift in time, i.e., the approximation has the same quality as the one given by a standard non-symplectic algorithm. Conversely, \cref{thm:RTSHamiltonian} proves that random step sizes do not spoil, under the assumptions specified above, the good long time properties of symplectic integrators with fixed step size.
\end{remark}
\begin{remark} As it can be noticed in the proof of \cref{lem:Remainder}, we introduce the assumption $p \geq 3/2$ in order to simplify the terms composing the remainder $S(\Delta, \eta)$. In case $1 \leq p < 3/2$, e.g. when the symplectic Euler method is employed ($q = 1$) and the natural scaling $p = q$ is chosen, the $\OO(h^q)$ approximation of the Hamiltonian still holds but with a slight reduction in the exponential terms appearing in the time span of validity.
\end{remark}
\begin{remark} Let us remark that in order for \eqref{eq:BoundS} to hold we implicitly assumed $t_n \geq 1$ to bound $\sqrt{t_n} \leq t_n$. If $t_n < 1$, we can bound every appearance of $t_n$ from \eqref{eq:BoundFirstTerm} to \eqref{eq:BoundS} as $t_n \leq 1$, and the desired result would still hold.
\end{remark}

\section{Bayesian inference}\label{sec:BayesianInference} It has been recently shown \cite{CGS17, CCC16, LST18} that probabilistic methods for ordinary and partial differential equations guarantee robust results (with respect to the numerical discretization error) in the context of Bayesian inverse problems. In this section, we briefly introduce a Bayesian inverse problem in the ODE setting and illustrate how the RTS-RK method can be employed in this framework. 

Let us consider a function $f_\theta \colon \R^d \to \R^d$ which depends on a real parameter $\theta \in \Theta$, where $\Theta$ is an open subset of $\R^n$ and the ODE
\begin{equation}\label{eq:ODEParam}
	y_\theta' = f_\theta(y), \quad y_\theta(0) = y_0 \in \R^d.
\end{equation}
In order to simplify the notation, we consider $y_0$ to be a fixed initial condition. In general, $y_0$ could depend itself on $\theta$. In the classical setting of numerical analysis, the main problem of interest is to determine the solution $y_\theta$ given the parameter $\theta$. The inverse problem we consider is instead to determine $\theta$ through observations of the solution $y_\theta$ (or quantities derived from it). In the Bayesian setting, the inverse problem is recast in terms of probability distributions, and the goal is to establish a probability measure on $\theta$, known as the posterior measure, given observed data and a probability measure, known as the prior, which captures all knowledge on the parameter available beforehand.

Let us denote by $z \in \R^m$ the observable and by $\mathcal G\colon \Theta \to \R^m$ the forward operator, which can be written as $\mathcal G = \mathcal O \circ \mathcal S$, where $\mathcal S$ is the solution operator and $\mathcal O$ is the observation operator. In this case, $\mathcal S \colon \R^n \to \mathcal C([0, T])$ is the operator mapping $\theta$ into the solution $y_\theta$, and $\mathcal O \colon \mathcal C([0, T]) \to \R^m$ maps the solution into the observable. Observations are then given by evaluations of the forward model corrupted by noise. In particular, we model noise as a Gaussian random variable $\epl \sim \mathcal{N}(0, \Sigma_\epl)$ independent of $\theta$, so that observations read
\begin{equation}
	z = \mathcal{G}(\theta) + \epl.
\end{equation}
Under these assumptions, the likelihood of the observations can be written as
\begin{equation}\label{eq:Likelihood}
	\pi(z \mid \theta) = e^{-V_z(\theta)},
\end{equation}
where the function $V_z\colon\Theta \to \R$, called the potential or negative log-likelihood, is given by
\begin{equation}\label{eq:Potential}
	V_z(\theta) = \frac{1}{2}\big(\mathcal{G}(\theta) - z\big)^\top \Sigma_\epl^{-1} \big(\mathcal{G}(\theta) - z\big).
\end{equation}

The second building block of Bayesian inverse problems is the prior distribution, which we denote by $\pi_0(\theta)$. The prior encodes all the knowledge on the parameter that is known before observations are provided. In the following, we adopt a common abuse of notation, confounding measures and their probability density function. 

Once the likelihood model and the prior distribution are established, it is possible to compute the posterior distribution $\pi(\theta \mid z)$ via Bayes' theorem, i.e.,
\begin{equation}
	\pi(\theta\mid z) = \frac{\pi(z\mid \theta)\pi_0(\theta)}{\mathcal Z(z)},
\end{equation}
where $\mathcal Z(z)$ is the normalising constant given by
\begin{equation}
	\mathcal Z(z) = \int_{\Theta} \pi(z\mid \theta) \pi_0(\theta) \, \dd \theta.
\end{equation}
Let us denote by $\mathcal{G}^h(\theta)$ the forward model where the solution operator is approximated by a Runge--Kutta method with time step $h$, and consequently with $V^h_z(\theta)$ and $\pi^h(z\mid\theta)$ the potential and the likelihood function obtained replacing $\mathcal{G}(\theta)$ with $\mathcal{G}^h(\theta)$. We can then define analogously the approximated posterior distribution $\pi^h(\theta \mid z)$ via Bayes' formula. In the following, we assume that the posteriors $\pi(\theta \mid z)$ and $\pi^h(\theta \mid z)$ are absolutely continuous with respect to the Lebesgue density. In \cite[Theorem 4.6]{Stu10}, Stuart proves that the posterior distribution $\pi^h(\theta\mid z)$ converges to $\pi(\theta\mid z)$ with respect to $h$ with the same rate as $V^h_z(\theta)$ converges to $V_z(\theta)$. There, convergence is shown with respect to the Hellinger distance for a Gaussian prior, which is defined for probability density functions which are absolutely continuous with respect to the Lebesgue density as
\begin{equation}
	\Hell\big(\pi^h(\theta\mid z), \pi(\theta\mid z)\big)^2 = \frac{1}{2}\int_{\Theta} \Big(\sqrt{\pi^h(\theta\mid z)} - \sqrt{\pi(\theta\mid z)}\Big)^2 \, \dd \theta.
\end{equation}
Hence, when there is no restriction in computational resources and it is possible to choose $h$ small, the approximated posterior distribution can be made arbitrarily close to the true posterior. The result is proved in \cite{Stu10} under the hypothesis of a Gaussian prior, but can be extended to a wider class of thin-tailed priors as done in \cite{DaS16} and to heavy-tailed priors as done in \cite{Sul17}.

In this work we consider the case when $h$ is fixed, and in particular we are interested in the case where the numerical error dominates the noise contribution. It has been shown via examples in \cite{CGS17, COS17} that in this small noise limit the approximated posterior distributions can be overly confident on the value of the parameter. In particular, the expectation of $\theta$ computed under the posterior distribution exhibits a bias with respect to the true value, which is not highlighted by the dispersion of the posterior itself. This undesirable phenomenon can be corrected by means of a probabilistic method, as the one presented by Conrad et al. in \cite{CGS17} or the RTS-RK method, to approximate the potential $V_z(\theta)$. Let us denote by $\xi \in \mathcal{X}$ the auxiliary random variable introduced by the probabilistic method. In the case of RTS-RK, we have $\xi = (H_0, H_1, \ldots, H_{N-1})^\top$ and $\mathcal{X} \subset \R_+^N$. The likelihood function, denoted as  $\pi^h_{\mathrm{pr}}(z\mid\theta)$ is then defined by
\begin{equation}\label{eq:LikelihoodProb}
	\pi^h_{\mathrm{prob}}(z\mid \theta) = \E^\xi e^{-V^{h,\xi}_{z}(\theta)}.
\end{equation}
where $V^{h, \xi}_{z}$ is the approximation of the potential function given by the probabilistic method. The corresponding posterior distribution $\pi^h_{\mathrm{pr}}$ is then defined by
\begin{equation}\label{eq:PosteriorProb}
	\pi^h_{\mathrm{prob}}(\theta\mid z) = \frac{\pi^h_{\mathrm{prob}}(z\mid \theta)\pi_0(\theta)}{\E^\xi \mathcal Z^{h, \xi}(z)},
\end{equation}
where the normalising constant is given by $\E^\xi \mathcal Z^{h, \xi}(z)$, where
\begin{equation}
	\mathcal Z^{h, \xi}(z) = \int_{\Theta} e^{-V^{h,\xi}}_{z}(\theta) \pi_0(\theta) \dd \theta.
\end{equation}
Modifying the posterior in this manner allows to obtain qualitatively better results, which account for the uncertainty introduced by the numerical solver. Moreover, this posterior distribution still converges to the true posterior for $h \to 0$ as proved in \cite{LST18}, where \eqref{eq:PosteriorProb} is called the marginal posterior.

In order to sample from the posteriors defined above we employ Markov chain Monte Carlo (MCMC) algorithms. In particular, due to the manner in which the probabilistic posterior \eqref{eq:PosteriorProb} is defined, the pseudo-marginal Metropolis--Hastings (PMMH) algorithm \cite{AnR09} is a suitable choice for sampling. We note that in case of a deterministic approximation of the forward model, the standard random walk Metropolis--Hastings can be employed.

\subsection{Analytical posteriors in a linear problem}\label{sec:AnalyticalPosterior} 

\begin{figure}[t!]
	\begin{center}
		\begin{tikzpicture}
		\draw (0,0) -- (9.2, 0) -- (9.2, 0.5) -- (0, 0.5) -- (0, 0);
		\draw[color=leg1,style=dotted,thick] (0.1, 0.25) -- (0.6, 0.25) node[right,color=black] {\small $\sigma=0.1$};
		\draw[color=leg2,style=dashdotted,thick] (2.1, 0.25) -- (2.6, 0.25) node[right,color=black] {\small$\sigma=0.05$};
		\draw[color=leg3,style=dashed,thick] (4.3, 0.25) -- (4.9, 0.25) node[right,color=black] {\small $\sigma=0.025$};			
		\draw[color=leg4,thick] (6.7, 0.25) -- (7.3, 0.25) node[right,color=black] {\small $\sigma=0.0125$};
		\end{tikzpicture}
		
		\vspace{0.2cm}
		\begin{tabular}{c@{\hskip 0.5in}c}
			\includegraphics[]{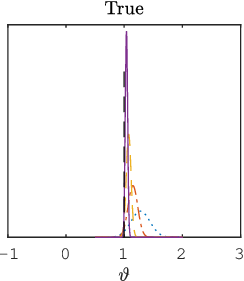} & \includegraphics[]{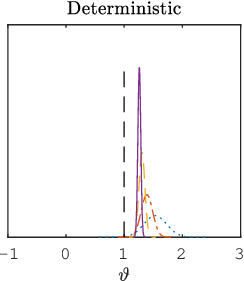} \\[10pt] 
			\includegraphics[]{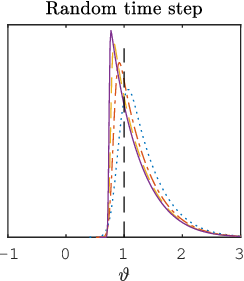}  & \includegraphics[]{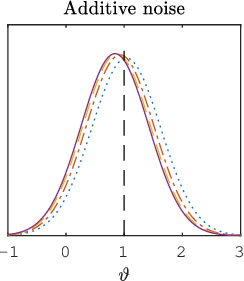}
		\end{tabular}	
	\end{center}
	\caption{Analytical posterior distributions in the linear case of \cref{sec:AnalyticalPosterior} for the true solution and its approximations with the deterministic explicit Euler method and the two probabilistic versions with additive noise \eqref{eq:AN-RK} and with random time steps \eqref{eq:RTS-RK}. In this case, $h = 0.5$ and the variance $\sigma^2$ of the observation error is reduced progressively. The true value of the initial condition $\theta^* = 1$ is shown with a vertical black dashed line.}
	\label{fig:AnalyticalPosterior}
\end{figure}

If the forward operator $\mathcal{G}$ is linear, the prior on the unknown parameter is Gaussian, and the negative log-likelihood is given by \eqref{eq:Potential}, then there is an explicit formula for the corresponding posterior distribution. Let us hence consider the following one dimensional ODE
\begin{equation}
	y'(t) = -y(t), \quad y(0) = \theta.
\end{equation}
Given $h > 0$, we consider the inferential problem of determining the true initial condition $\theta^*$ from a single observation $z = \phi_h(\theta^*) + \epl$, where $\phi_h(\theta^*) = \theta^*e^{-h}$ is the true solution at time $t =  h$ and $\epl \sim \mathcal N (0, \sigma^2)$ is a source of noise. In this case, the parameter space is $\Theta = \R$ and the forward operator $\mathcal G$ is defined by $\mathcal G\colon \R \to \R$, $\mathcal G \colon \theta \mapsto \theta e^{-h}$. In the following, we verify heuristically the convergence of the posterior distributions obtained with deterministic and probabilistic integrators with respect to a vanishing noise scale. If a Gaussian prior $\pi_0 = \mathcal N(0, 1)$ is given for $\theta$, the true posterior distribution is computable analytically and is given by
\begin{equation}\label{eq:ExPosteriorEx}
	\pi(\theta \mid z) = \mathcal N\Big(\theta; \frac{ze^{-h}}{\sigma^2 + e^{-2h}}, \frac{\sigma^2}{\sigma^2 + e^{-2h}}\Big),
\end{equation}
where $\mathcal N(x; \mu, \alpha^2)$ is the density of a Gaussian random variable of mean $\mu$ and variance $\alpha^2$ evaluated in $x$. Consistently, if $\sigma^2 \to 0$, we have that $z \to \theta^*e^{-h}$ and therefore $\pi(\theta \mid z) \to \delta_{\theta^*}$. 

If we approximate $\phi_h(\theta)$ for a given initial condition $\theta$ with a single step of the explicit Euler method (i.e., with step size $h$), we get $\Psi_h(\theta) = (1 - h)\theta$. Computing the posterior distribution obtained with this approximation leads to 
\begin{equation}\label{eq:ExPosteriorRK}
\pi^h(\theta \mid z) = \mathcal N\Big(\theta; \frac{(1 - h)z}{\sigma^2 + (1- h)^2}, \frac{\sigma^2}{\sigma^2 + (1 - h)^2}\Big).
\end{equation}
In the limit of $\sigma^2 \to 0$, we get in this case that the posterior distribution tends to $\pi^h(\theta \mid z) \to \delta_{\bar \theta}$, where $\bar \theta = e^{-h}\theta^* / (1 - h)$. The posterior distribution is hence tending to a biased Dirac delta with respect to the true value.

Let us consider the additive noise method \eqref{eq:AN-RK} applied to the explicit Euler method, i.e., the random approximation $y(h) \approx Y_1$, where $Y_1 = (1 - h)\theta + \xi$ and where $\xi \sim \mathcal N(0, h^3)$, so that the method converges consistently with the deterministic method. In this case, the posterior distribution that we denote by $\pi^h_{\mathrm{prob, AN}}$ is given by
\begin{equation}
\pi^h_{\mathrm{prob, AN}}(\theta \mid z) = \mathcal N\Big(\theta; \frac{(1 - h)z}{\tilde \sigma^2 + (1- h)^2}, \frac{\tilde \sigma^2}{\tilde \sigma^2 + (1 - h)^2}\Big).
\end{equation}
where $\tilde \sigma^2 = \sigma^2 + h^3$. Hence, taking the limit $\sigma^2 \to 0$ gives
\begin{equation}\label{eq:ExPosteriorAN}
\pi^h_{\mathrm{prob, AN}}(\theta \mid z) \to \mathcal N\Big(\theta; \frac{(1 - h)e^{-h}\theta^*}{h^3 + (1 - h)^2}, \frac{h^3}{h^3 + (1 - h)^2}\Big),
\end{equation}
which shows that while the asymptotic mean is still biased with respect to the true value, the uncertainty in the forward model is reflected by a positive variance.
Let us now consider the random time step explicit Euler with step size distribution $H \sim \mathcal U(h - h^{p+1/2}, h + h^{p+1/2})$. In this case, the forward model is given by
\begin{equation}
Y_1 = \theta - H \theta = (1 - h)\theta + U \theta, \quad U \sim \mathcal U (-h^{p+1/2}, h^{p+1/2}).
\end{equation}
Hence, disregarding all multiplicative constants that are independent of $\theta$ and setting $p = q = 1$, we get the posterior
\begin{equation}\label{eq:ExPosteriorRTS}
\pi^h_{\mathrm{prob, RTS}}(\theta \mid z) \propto \exp\Big(-\frac{\theta^2}{2}\Big) \frac{1}{\theta} \Big(\Phi\Big(\frac{((1 - h) + h^{3/2})\theta - z}{\sigma}\Big) - \Phi\Big(\frac{((1 - h) - h^{3/2})\theta - z}{\sigma}\Big) \Big),
\end{equation}
where $\Phi$ denotes the cumulative distribution function of a standard Gaussian random variable. Since we require in  \cref{as:hStrong}.\ref{as:hStrong_Pos} that $H > 0$ almost surely, the time step $H$ cannot be Gaussian and the closed-form expression of the posterior is not as neatly defined as in the additive noise case. In the limit for $\sigma \to 0$, we get the limiting distribution 
\begin{equation}
\pi^h_{\mathrm{prob, RTS}}(\theta \mid z) \propto \exp\Big(-\frac{\theta^2}{2}\Big)\frac{1}{\theta} \chi_{\{y_{\min} \leq \theta \leq y_{\max}\}},
\end{equation}
where $y_{\min}$ and $y_{\max}$ are given by
\begin{equation}
y_{\min} = \frac{e^{-h}\theta^*}{((1 - h) + h^{3/2})}, \quad y_{\max} = \frac{e^{-h}\theta^*}{((1 - h) - h^{3/2})}.
\end{equation}
It is hence possible to remark that for the RTS-RK method the variance of the posterior distribution is not collapsing to zero for $\sigma \to 0$ as in the deterministic case. 

We fix $h = 0.5$ and consider $\sigma = \{0.1, 0.05, 0.025, 0.0125\}$, thus generating four observational noises $\eta_i$ as $\eta_i = \sigma_i Z$ for a random variable $Z \sim \mathcal{N}(0, 1)$. In \cref{fig:AnalyticalPosterior} we show the posteriors \eqref{eq:ExPosteriorEx}, \eqref{eq:ExPosteriorRK}, \eqref{eq:ExPosteriorAN} and \eqref{eq:ExPosteriorRTS}, which confirm our claim, i.e., that probabilistic methods take into account the variability in the forward model caused by the numerical approximation and transfer it to the posterior belief.

\section{Numerical experiments}\label{sec:NumericalExperiments}

In this section, we present a series of numerical experiments that illustrate the versatility and usefulness of our new random time stepping method. These experiments also corroborate the theoretical results presented in the previous sections.

\subsection{Weak order of convergence}\label{sec:NumericalExperimentsWeak}

\begin{figure}[t!]
	\centering
	\begin{tabular}{c@{\hspace{0.3cm}}c}
		\includegraphics[]{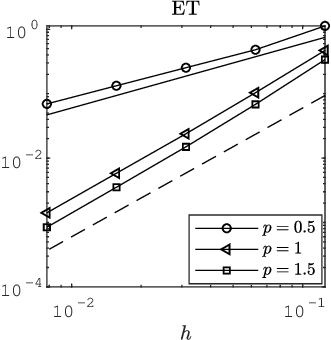} & \includegraphics[]{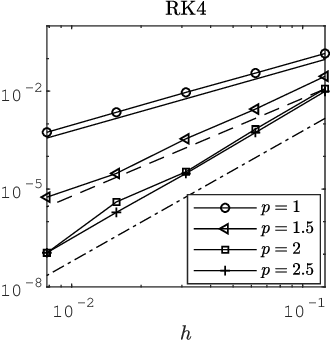} \\
	\end{tabular}
	\caption{Weak order of convergence for the random time-stepping explicit trapezoidal (ET) and fourth-order Runge--Kutta (RK4) as a function of the value of $p$ of \cref{as:hStrong}. In the left figure, reference slopes $1$ and $2$ are displayed (solid and dashed lines), while in the right figure reference slopes $2$, $3$ and $4$ are displayed (solid, dashed and dash-dotted lines).}
	\label{fig:Weak}
\end{figure}

In order to verify the result predicted in \cref{thm:StrongOrder}, we consider the FitzHugh--Nagumo equation, which is defined as
\begin{equation}\label{eq:FitzNag}
\begin{aligned}
y_1' &= c\big(y_1 - \frac{y_1^3}{3} + y_2\big), && y_1(0) = -1, \\
y_2' &= -\frac{1}{c}(y_1 - a + by_2), && y_2(0) = 1,
\end{aligned}
\end{equation}
where $a, b, c$ are real parameters with values $a = 0.2$, $b = 0.2$, $c = 3$. We integrate the equation from time $t_0 = 0$ to final time $T = 1$. The reference solution is generated with a high-order method on a fine time scale. The deterministic integrators we choose in this experiment are the explicit trapezoidal rule and the classic fourth-order Runge--Kutta method. The random steps are uniform as in \cref{ex:uniformH}. We vary their mean in the range $h_i = 0.125\cdot 2^{-i}$ with $i = 0, 1, \ldots, 4$, and we vary the value of $p$ in \cref{as:hStrong} in order to verify the theoretical result of \cref{thm:weakOrder}. In particular, we consider $p \in \{0.5, 1, 1.5\}$ for the explicit trapezoidal rule and $p \in \{1, 1.5, 2, 2.5\}$ for the classic fourth order Runge--Kutta method. The function $\Phi\colon\R^d\to\R$ of the solution we consider is defined as $\Phi(x) \defeq x^\top x$. Finally, we consider $10^6$ trajectories of the numerical solution in order to approximate the expectation with a Monte Carlo sum. Results (\cref{fig:Weak}) show that the order of convergence predicted theoretically is confirmed by numerical experiments. 

\subsection{Mean square order of convergence}\label{sec:NumericalExperimentsStrong}

\begin{figure}[t!]
	\centering
	\begin{tabular}{c@{\hspace{0.3cm}}c}
		\includegraphics[]{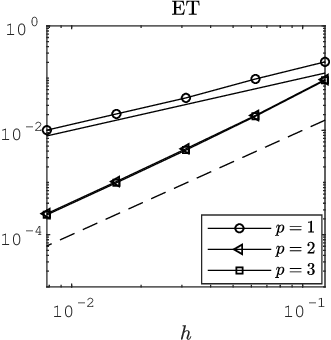} & \includegraphics[]{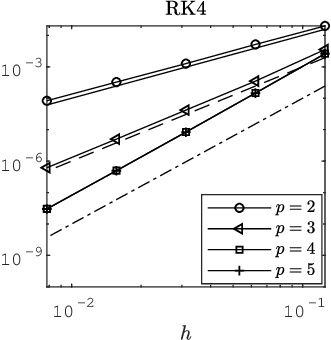} \\
	\end{tabular}
	\caption{Mean square order of convergence for the random time-stepping explicit trapezoidal (ET) and fourth-order Runge--Kutta (RK4) as a function of the value of $p$ of \cref{as:hStrong}. In the left figure, reference slopes $1$ and $2$ are displayed (solid and dashed lines), while in the right figure reference slopes $2$, $3$ and $4$ are displayed (solid, dashed and dash-dotted lines).}
	\label{fig:MeanSquare}
\end{figure}

We now verify the weak order of convergence predicted in \cref{thm:weakOrder}. For this experiment we consider the ODE \eqref{eq:FitzNag} as well, with the same time scale $T$ and parameters as in \cref{sec:NumericalExperimentsWeak}. The reference solution at final time is generated in this case as well with a high-order method on a fine time scale. We consider as deterministic solvers the explicit trapezoidal rule and the classic fourth order Runge--Kutta method, which verify \cref{as:PsiStrong} with $q = 2$ and $q = 4$ respectively. Moreover, we consider uniform random time steps as in \cref{ex:uniformH}, where we vary the value of $p$ in \cref{as:hStrong} in order to verify the order of convergence predicted in \cref{thm:StrongOrder}. In particular, we consider $p \in \{1, 2, 3\}$ for the explicit trapezoidal rule and $p \in \{2, 3, 4, 5\}$ for the classic fourth order Runge--Kutta method. We vary the mean time step $h$ taken by the random time steps $H_n$ in the range $h_i = 0.125\cdot 2^{-i}$, with $i = 0, 1, \ldots, 4$. Then, we simulate $10^3$ realizations of the numerical solution $Y_{N_i}$, with $N_i = T / h_i$ for $i = 0, 1, \ldots, 4$, and compute the approximate mean square order of convergence for each value of $h$ with a Monte Carlo mean. Results (\cref{fig:MeanSquare}) show that the orders predicted theoretically by \cref{thm:StrongOrder} are confirmed numerically. 

\begin{figure}[t!]
	\centering
	\begin{tabular}{c@{\hspace{0.3cm}}c}
		\includegraphics[]{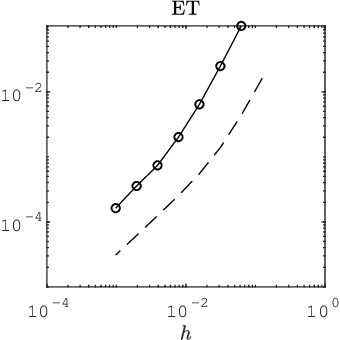} & \includegraphics[]{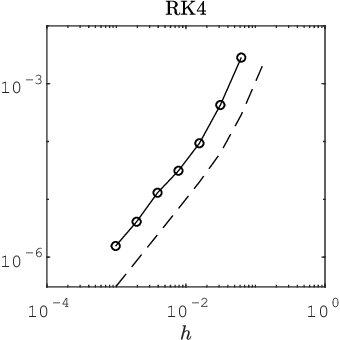} \\
	\end{tabular}
	\caption{Convergence of the {square root of the} MSE of the Monte Carlo estimator for the random time-stepping explicit trapezoidal (ET) (left figure) and fourth-order Runge--Kutta (RK4) (right figure) {with respect to the time step $h$}. The dashed line corresponds to the order predicted in \cref{thm:MSEMonteCarlo} with $M = 10^3$ for ET and $M = 10^4$ for RK4.}
	\label{fig:MonteCarlo}
\end{figure}

\begin{figure}[t!]
	\centering
	\begin{tabular}{c@{\hspace{0.3cm}}c}
		\includegraphics[]{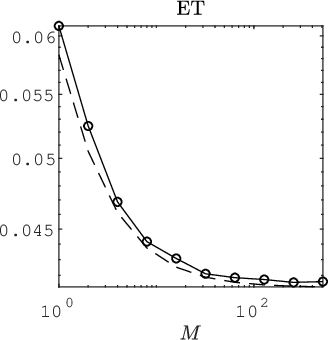} & \includegraphics[]{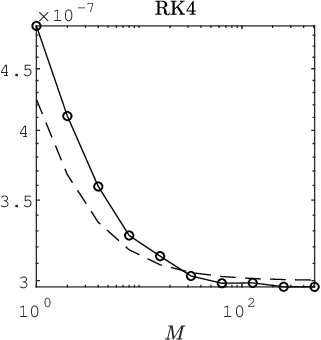} \\
	\end{tabular}
	\caption{{Convergence of the square root of the MSE of the Monte Carlo estimator for the random time-stepping explicit trapezoidal (ET) (left figure) and fourth-order Runge--Kutta (RK4) (right figure) with respect to the number of trajectories $M$. The dashed line corresponds to the order predicted in \cref{thm:MSEMonteCarlo} with $h = 0.05$ for ET and $h = 0.01$ for RK4.}}
	\label{fig:MonteCarlo_M}
\end{figure}

\subsection{Mean-square convergence of Monte Carlo estimators} 
We shall now verify numerically the validity of \cref{thm:MSEMonteCarlo}. We consider the ODE \eqref{eq:FitzNag}, with final time $T = 1$ and the same parameters as above. In this case as well, we consider the explicit trapezoidal rule and the fourth-order explicit Runge--Kutta method with uniform random time steps having mean $h_i = 0.125\cdot 2^{-i}$ with $i = 0, 1, \ldots, 7$. For the explicit trapezoidal rule, we fix $M = 10^3$ and $p = 1$, so that for bigger values of $h$ the first term in the bound presented in \cref{thm:MSEMonteCarlo} dominates, while in the regime of small $h$, the higher order of the first term makes the second term larger in magnitude. This behaviour results in the change of slope in the convergence plot which can be observed in \cref{fig:MonteCarlo}, both in the theoretical estimate and in the numerical results. We perform the same experiment using the fourth-order explicit Runge--Kutta method, fixing $M = 10^4$ and $p = 1.5$, thus obtaining a numerical confirmation of the theoretical result.

{As a second experiment, we consider the same setup as above but wish to verify the dependence of the MSE on the number of samples $M$, which we vary as $M = 2^i$, with $i = 0, 1, \ldots, 9$. For the explicit trapezoidal rule, we consider $p = q = 2$, which is the optimal choice for the intrinsic variability of the RTS-RK method. Moreover, we fix $h = 0.05$. In this case, the bound \eqref{eq:MSEBound} reduces to
\begin{equation}
	\mathrm{MSE}(\widehat Z_{N, M}) \leq Ch^{2q} \Big(1 + \frac{1}{M}\Big).
\end{equation}
In \cref{fig:MonteCarlo_M} we show that the convergence of the MSE depends on $M$ as predicted by the theoretical bound. We repeat the same experiment using the fourth order explicit Runge--Kutta method, for which we take $h = 0.01$ and $p = q = 4$, thus confirming numerically our theoretical result.}

\subsection{Robustness} In this numerical experiment we verify the robustness of RTS-RK when applied to chemical reactions. Let us consider the Peroxide-Oxide chemical reaction, which is macroscopically defined by the following balance equation
\begin{equation}
	\mathrm{O}_2 + 2\mathrm{NADH} + 2\mathrm{H}^+ \to 2\mathrm{H}_2\mathrm{O} + 2\mathrm{NAD}^+,
\end{equation}
where NADH and NAD$^+$ are the oxidized and reduced form of the nicotinamide adenine dinucleotide (NAD) respectively. This reaction has to be catalyzed by an enzyme to take place, which reacts with the reagents to create intermediate products of the reaction. A successful model \cite{Ols83} to describe the time-evolution of the chemical system is the following
\begin{equation}
\begin{aligned}
	\mathrm{B} + \mathrm{X} &\rightarrowtext{k_1} 2 \mathrm{X}, 
	&&2\mathrm{X} \rightarrowtext{k_2} 2\mathrm{Y}, 
	&&\mathrm{A} + \mathrm{B} + \mathrm{Y} \rightarrowtext{k_3} 3 \mathrm{X}, \\
	\mathrm{X} &\rightarrowtext{k_4} \mathrm{P}, 
	&&\mathrm{Y} \rightarrowtext{k_5} \mathrm{Q}, 
	&&\mathrm{X_0} \rightarrowtext{k_6} \mathrm{X}, \\
	\mathrm{A_0} &\leftrightarrowtext{k_7} \mathrm{A}, 
	&&\mathrm{B_0} \rightarrowtext{k_8} \mathrm{B}.
\end{aligned}
\end{equation}
Here, A and B are respectively $[\mathrm{O}_2]$ and $[\mathrm{NADH}]$, P, Q are the products and X, Y are intermediate results of the reaction process. It is therefore possible to model the time evolution of the reaction with the following system of nonlinear ODEs 
\begin{equation}\label{eq:PeroxOx}
\begin{aligned}
	\mathrm{A}' &= k_7  (\mathrm{A}_0 - \mathrm{A}) - k_3  \mathrm{A}\mathrm{B}\mathrm{Y}, &&\mathrm{A}(0) = 6, \\
	\mathrm{B}' &= k_8\mathrm{B}_0 - k_1  \mathrm{B}\mathrm{X} - k_3  \mathrm{A}\mathrm{B}\mathrm{Y}, &&\mathrm{B}(0) = 58, \\
	\mathrm{X}' &= k_1  \mathrm{B}\mathrm{X} - 2  k_2  \mathrm{X}^2 + 3  k_3 \mathrm{A}\mathrm{B}\mathrm{Y} - k_4  \mathrm{X} + k_6\mathrm{X}_0,&& \mathrm{X}(0) = 0, \\
	\mathrm{Y}' &= 2  k_2  \mathrm{X}^2 - k_5  \mathrm{Y} - k_3  \mathrm{A}\mathrm{B}\mathrm{Y}, && Y(0) = 0,
\end{aligned}
\end{equation}
where $\mathrm{A}_0 = 8$, $\mathrm{B}_0 = 1$, $\mathrm{X}_0 = 1$ and the real parameters $k_i$, $i = 1, \ldots, 8$ representing the reaction rates take values
\begin{equation}
\begin{aligned}
k_1 &= 0.35, &&k_2 = 250, &&k_3 = 0.035, &&k_4 = 20,\\
k_5 &= 5.35, &&k_6 = 10^{-5}, &&k_7 = 0.1, &&k_8 = 0.825.
\end{aligned}
\end{equation}            
It has been shown \cite{Ols83} that for these values of the parameters the system exhibits a chaotic behavior. In particular, at long time the trajectories lie in a strange attractor, and the system shows a strong sensitivity to perturbations on the initial condition. 

Since the components of the solution represent the concentration of chemicals, we require the numerical solution to be positive. Apart from physical considerations, numerically we observe that if one of the components takes negative values, the solution shows strong instabilities. For the RTS-RK method, the distribution of the random time steps can be selected so that the probability of obtaining a negative solution is zero, see e.g. \cref{ex:uniformH}. In contrast, for the additive noise method we can have disruptive effects even for $h$ small if the solution has a small magnitude, as the probability for negative populations will never be zero. Hence, in this case employing the additive noise method likely produces instabilities regardless of the chosen time step.

\begin{figure}
	\begin{center} 
		\begin{tabular}{c}
			\includegraphics[]{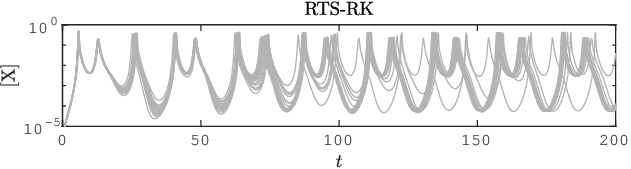}\\[10pt]
			\includegraphics[]{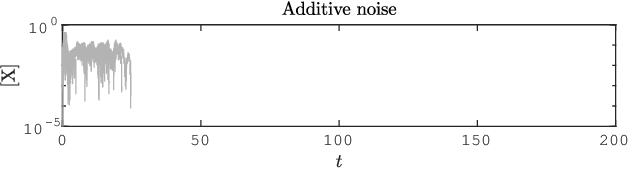}
		\end{tabular}
	\end{center}
	\caption{Fifty trajectories of the numerical value of the concentration of the $\mathrm X$ species for the random time-stepping and additive noise methods (above and below respectively).}
	\label{fig:OxPeroxTraj}
\end{figure}

Let us apply the additive noise method \eqref{eq:AN-RK} and the random time-stepping scheme \eqref{eq:RTS-RK} to equation \eqref{eq:PeroxOx}. We choose $h = 0.05$ as the mean of uniformly distributed time steps for \eqref{eq:RTS-RK} and as the time step for \eqref{eq:AN-RK}, while we employ the Runge--Kutta--Chebyshev method (RKC) \cite{HoS80} as deterministic integrator. Since RKC has order 1, we fix $p = q = 1$. As the problem is stiff, stabilized methods prevent a step size restriction while remaining explicit. We note that the RKC method is a stabilized numerical integrator of first order and that higher order explicit stabilized methods such as ROCK2 or ROCK4 \cite{AbM01, Abd02} could also be used as deterministic solvers for the RTS-RK method. It can be seen in \cref{fig:OxPeroxTraj} that the RTS-RK method conserves the positivity of the numerical solution while capturing the chaotic nature of the chemical reaction. In contrast, the additive noise scheme produces negative values, thus showing strong instabilities in the long-time behavior. In particular, all the numerical trajectories turn negative or diverge before approximately $t = 25$, which is the reason why after this time they are not displayed in \cref{fig:OxPeroxTraj}.

\subsection{Conservation of quadratic first integrals} 

\begin{figure}[t]
	\begin{center}
		\begin{tabular}{cccc}
			\includegraphics[]{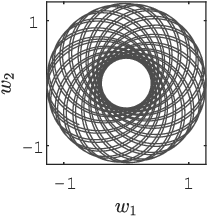} & \includegraphics[]{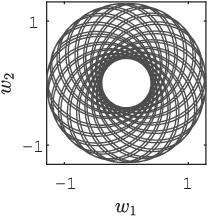} \includegraphics[]{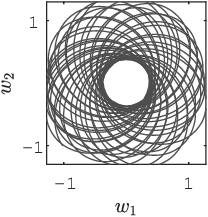} & \includegraphics[]{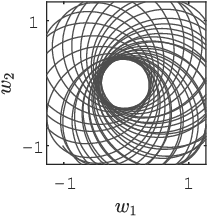} \\
		\end{tabular}\\
	\vspace{0.4cm}
	\includegraphics[]{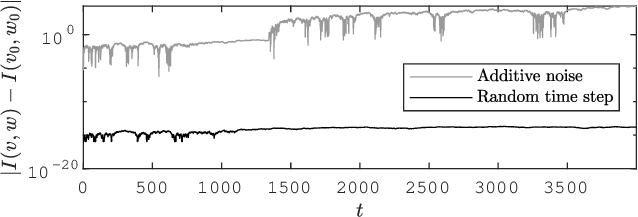}
	\end{center}
	\caption{Trajectories of \eqref{eq:KeplerPert} given by the RTS-RK method \eqref{eq:RTS-RK} for $0 \leq t \leq 200$ and $3800 \leq t \leq 4000$ (first and second figures), and by the additive noise method \eqref{eq:AN-RK} for $0 \leq t \leq 200$ and $200 \leq t \leq 400$ (third and fourth figures). Error on the angular momentum $I$ defined in \eqref{eq:AngularMomentum} for $0 \leq t \leq 4000$ given by the two methods.}
	\label{fig:Kepler}
\end{figure}

A simple model for the two-body problem in celestial mechanics is the Kepler system with a perturbation, which reads
\begin{equation}\label{eq:KeplerPert}
\begin{aligned}
	w_1' &= v_1, && v_1' = -\frac{w_1}{\norm{q}^3} - \frac{\delta w_1}{\norm{q}^5}, \\
	w_2' &= v_2, && v_2' = -\frac{w_2}{\norm{q}^3} - \frac{\delta w_2}{\norm{q}^5},
\end{aligned}
\end{equation}
where $v_1$, $v_2$ are the two components of the velocity and $w_1$, $w_2$ are the two components of the position. We set the perturbation parameter $\delta$ to be equal to 0.015 and the initial condition to be
\begin{equation}
	w_1(0) = 1 - e,\quad w_2(0) = 0, \quad v_1(0) = 0, \quad v_2(0) = \sqrt{(1 + e)/(1 − e)},
\end{equation}
where $e = 0.6$ is the eccentricity. It is well-known that this equation has the Hamiltonian and the angular momentum as quadratic first integrals. In particular, we focus here on the angular momentum, which reads
\begin{equation}\label{eq:AngularMomentum}
	I(v, w) = w_1v_2 - w_2v_1.
\end{equation}
We consider the simplest Gauss collocation method, namely the implicit midpoint rule, as the deterministic Runge--Kutta method. It is known that Gauss collocation methods conserve quadratic first integrals. According to \cref{thm:PolyInvariants}, we expect therefore that the random time-stepping method \eqref{eq:RTS-RK} implemented with $\Psi_h$ given by the implicit midpoint rule also conserves quadratic first integrals. We {integrate} \eqref{eq:KeplerPert} with uniformly distributed random time steps with mean $h = 0.01$ from time $t = 0$ to time $t = 4000$ which corresponds to approximately $636$ revolutions of the system (long-time behavior). Since the implicit midpoint rule is of order $q = 2$, we choose $p = 2$ for the RTS-RK method. Moreover, we consider the additive noise method \eqref{eq:AN-RK} with $h = 0.01$, expecting that the first integral will not be conserved. We observe in \cref{fig:Kepler} that the method \eqref{eq:RTS-RK} conserves the angular momentum, while for the method \eqref{eq:AN-RK} the approximate conservation of the quadratic first integral shown in \eqref{eq:BiasQuadraticAddNoise} is lost when integrating \eqref{eq:KeplerPert} over long time.

\subsection{Conservation of Hamiltonians} Let us consider the pendulum problem, which is given by the Hamiltonian $Q \colon \R^2 \to \R$ defined by

\begin{figure}[th!]
	\begin{center}
	\begin{subfigure}{\linewidth}
	\begin{subfigure}{0.7\linewidth}
	\includegraphics[]{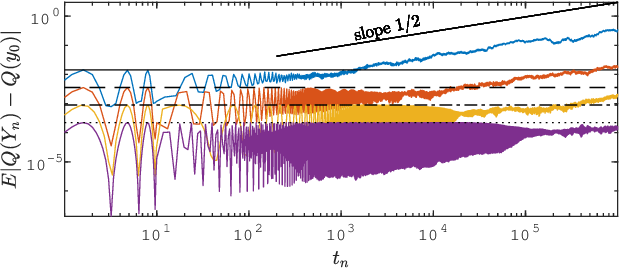} 
	\end{subfigure}
	\begin{subfigure}{0.2\linewidth}
	\begin{tikzpicture}
		\draw (0,0.5) -- (2.3, 0.5) -- (2.3, -2) -- (0, -2) -- (0, 0.5);
		\node at (0, 0.25) [right] {\small \textbf{(a)}};
		\draw[thick] (0.1, -0.25) -- (0.6, -0.25) node[right,color=black] {\small $h=0.2$};
		\draw[style=dashed,thick] (0.1, -0.75) -- (0.6, -0.75) node[right,color=black] {\small$h = 0.1$};
		\draw[style=dashdotted,thick] (0.1, -1.25) -- (0.6, -1.25) node[right,color=black] {\small $h=0.05$};			
		\draw[style=dotted,thick] (0.1, -1.75) -- (0.6, -1.75) node[right,color=black] {\small $h=0.025$};
	\end{tikzpicture}
	\end{subfigure}
	\end{subfigure}
	\begin{subfigure}{\linewidth}
	\begin{subfigure}{0.7\linewidth}
		\hspace{-0.1cm}\includegraphics[]{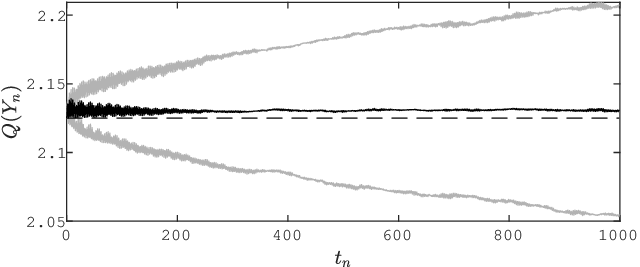} 
	\end{subfigure}
	\begin{subfigure}{0.2\linewidth}
	\begin{tikzpicture}
		\draw (0,0.5) -- (3.6, 0.5) -- (3.6, -2) -- (0, -2) -- (0, 0.5);
		\node at (0, 0.25) [right] {\small \textbf{(b1)}};
		\node at (0, -0.25) [right] {\small $h=0.2$};
		\draw[] (0.1, -0.75) -- (0.6, -0.75) node[right,color=black] {\small Mean};
		\draw[color=lightgray] (0.1, -1.25) -- (0.6, -1.25) node[right,color=black] {\small Confidence interval};
		\draw[style=dashed,thick] (0.1, -1.75) -- (0.6, -1.75) node[right,color=black] {\small $Q(y_0)$};
	\end{tikzpicture}
	\end{subfigure}	
	\end{subfigure}
	\begin{subfigure}{\linewidth}
	\begin{subfigure}{0.7\linewidth}
		\hspace{-0.1cm}\includegraphics[]{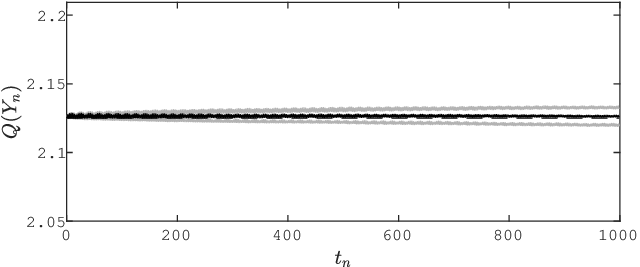} 
	\end{subfigure}
	\begin{subfigure}{0.2\linewidth}
		\begin{tikzpicture}
		\draw (0,0.5) -- (3.6, 0.5) -- (3.6, -2) -- (0, -2) -- (0, 0.5);
		\node at (0, 0.25) [right] {\small \textbf{(b2)}};
		\node at (0, -0.25) [right] {\small $h=0.1$};
		\draw[] (0.1, -0.75) -- (0.6, -0.75) node[right,color=black] {\small Mean};
		\draw[color=lightgray] (0.1, -1.25) -- (0.6, -1.25) node[right,color=black] {\small Confidence interval};
		\draw[style=dashed,thick] (0.1, -1.75) -- (0.6, -1.75) node[right,color=black] {\small $Q(y_0)$};
		\end{tikzpicture}
	\end{subfigure}
	\end{subfigure}
	\end{center}
	\caption{\textbf{(a)}: Time evolution of the mean error for the pendulum problem and different values of the time step $h$. The black lines represent the theoretical estimate given by \cref{thm:RTSHamiltonian}, while the colored lines represent the experimental results. The mean was computed {by} averaging 20 realisations of the numerical solution. {\textbf{(b1)} and \textbf{(b2)}: Time evolution of the mean  Hamiltonian for two different values of the time step. The mean Hamiltonian is depicted together with an approximate confidence interval, whose width is proportional to the standard deviation of the Hamiltonian over 200 trajectories.}}
	\label{fig:MeanTime}
\end{figure}

\begin{equation}
Q(v, w) = \frac{v^2}{2} - \cos w,
\end{equation}
where $y = (v, w)^\top \in \R^2$. We wish to study the validity of \cref{thm:RTSHamiltonian}, i.e., show that the mean error on the Hamiltonian is of order $\OO(h^q)$ for time spans of polynomial length and then it grows proportionally to the square root of time. We consider the initial condition $(v_0, w_0) = (1.5, -\pi)$ and integrate the equation employing RTS-RK based on the implicit midpoint method ($q = 2$) choosing $p = q$, which is the optimal scaling of the noise. We choose uniform time steps, vary their mean $h \in \{0.2, 0.1, 0.05, 0.025\}$, integrate the dynamical system up to the final time $T = 10^6$ and study the time evolution of the mean numerical error on the Hamiltonian $Q$. Results are shown in \cref{fig:MeanTime}, where it is possible to notice that the error is bounded by $\OO(h^q)$ (horizontal black lines) for long time spans. After this stationary phase, the error on the Hamiltonian appears to grow as the square root of time. The oscillations of the error which are shown in \cref{fig:MeanTime} are present even when integrating the pendulum system with a deterministic symplectic scheme. {Moreover, considering $T = 10^3$, the time step $h \in \{0.2, 0.1\}$ and keeping all other parameters as above, we compute the mean Hamiltonian and represent it in \cref{fig:MeanTime} together with an approximate confidence interval. We arbitrarily compute the confidence interval as $(\E Q(Y_n) - 2 \mathrm{Var}Q(Y_n)^{1/2}, \E Q(Y_n) + 2 \mathrm{Var}Q(Y_n)^{1/2})$, and we employ it to show the path-wise variability of the value of the Hamiltonian. As expected, the variability decreases dramatically with respect to the time step $h$.}

\subsection{Bayesian inference}\label{sec:BayesianInferenceEx}
\begin{figure}	
	\begin{center}
		\begin{tikzpicture}
		\draw (0,0) -- (9.2, 0) -- (9.2, 0.5) -- (0, 0.5) -- (0, 0);
		\draw[color=leg21,thick] (0.1, 0.25) -- (0.6, 0.25) node[right,color=black] {\small $h=0.2$};
		\draw[color=leg22,thick] (2.2, 0.25) -- (2.7, 0.25) node[right,color=black] {\small $h=0.1$};
		\draw[color=leg23,thick] (4.3, 0.25) -- (4.9, 0.25) node[right,color=black] {\small $h=0.05$};			
		\draw[color=leg24,thick] (6.7, 0.25) -- (7.3, 0.25) node[right,color=black] {\small $h=0.025$};
		\end{tikzpicture}
		
		\vspace{0.3cm}
		\begin{tabular}{ccc}
			\hspace{-0.32cm}\includegraphics[]{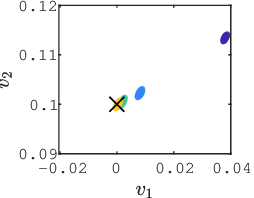}  & \hspace{-0.32cm}\includegraphics[]{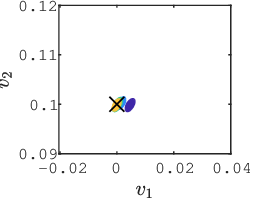}  & \hspace{-0.32cm}\includegraphics[]{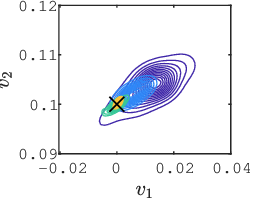} \\ 
			\includegraphics[]{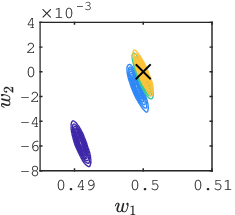} & \includegraphics[]{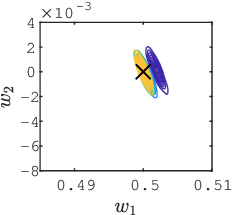} & \includegraphics[]{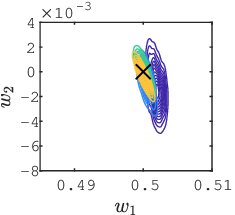} \\
		\end{tabular}
	\end{center}
	\caption{Posterior distributions for the initial position and velocity of the Hénon-Heiles system with different values of $h = \{0.2, 0.1, 0.05, 0.025\}$. First row: initial velocity $v_0$. Second row: initial position $w_0$. First column: deterministic Heun's method. Second column: deterministic Störmer--Verlet scheme. Third column: RTS-RK Störmer--Verlet ($p=2$).}
	\label{fig:Bayes}
\end{figure}

For the last numerical experiment we consider the Hénon--Heiles equation, a Hamiltonian system with energy $Q\colon \R^4 \to \R$ defined by
\begin{equation}\label{eq:HamHH}
	Q(v, w) = \frac{1}{2}\norm{v}^2 + \frac{1}{2}\norm{w}^2 + w_1^2w_2 - \frac{1}{3}w_2^3,
\end{equation} 
where $v, w \in \R^2$ are the velocity and position respectively and where we denote by $y = (v, w)^\top \in \R^4$ the solution. We consider an initial condition such that $Q(y_0) = 0.13$, for which the system exhibits a chaotic behaviour \cite{HeH64}. In the spirit of \cref{sec:BayesianInference}, we are interested in recovering the true value of the initial condition $y_0$ through a single observation $y_{\mathrm{obs}}$ of the solution $(v, w)$ at a fixed time $t_{\mathrm{obs}} = 10$. The exact forward operator $\mathcal{G}$ is therefore defined as $\mathcal{G}(y_0) = \phi_{t_{\mathrm{obs}}}(y_0)$. Noise is then set to be a Gaussian random variable $\epl \sim \mathcal{N}(0, \sigma_\epl^2 I)$, where $\sigma_\epl = 5 \cdot 10^{-4}$, and we fix a standard Gaussian prior on the initial condition, i.e., $\pi_0 = \mathcal N(0, I)$, so that the likelihood is given by \eqref{eq:Likelihood}. We choose the observational noise to have a small variance (i.e., of order $\mathcal O(10^{-8})$) as in this case classical solvers present the misleading overconfident behaviour explained in \cref{sec:BayesianInference}.

Since the equation is Hamiltonian, we choose to employ a classical second-order ($q = 2$) symplectic method, the Störmer--Verlet scheme \cite{Sto07, Ver67, HLW06}, for which one step is defined in the general case as
\begin{equation}
\begin{aligned}
v_{n+1/2} &= v_n - \frac{h}{2} \nabla_w Q(v_n, w_n), \\
w_{n+1} &= w_n + \frac{h}{2} \big(\nabla_v Q(v_{n+1/2}, w_n) + \nabla_v Q(v_{n+1/2}, w_{n+1})\big),\\
v_{n+1} &= v_{n+1/2} - \frac{h}{2} \nabla_w Q(v_{n+1/2}, w_{n+1}).
\end{aligned}
\end{equation}
As the Hamiltonian $Q$ given by \eqref{eq:HamHH} is separable, i.e., $Q(v, w) = Q_1(v) + Q_2(w)$, where $Q_1, Q_2 \colon \R^2 \to \R$, the Störmer--Verlet scheme simplifies to
\begin{equation}
\begin{aligned}
v_{n+1/2} &= v_n - \frac{h}{2} \nabla_w Q_2(w_n), \\
w_{n+1} &= w_n + h \nabla_v Q_1(v_{n+1/2}),\\
v_{n+1} &= v_{n+1/2} - \frac{h}{2} \nabla_w Q_2(w_{n+1}).
\end{aligned}
\end{equation}
Hence, in the separable case the Störmer--Verlet scheme is explicit and the evaluation of the flow consists only of three evaluations of the derivatives of $Q$. We then employ this method both with a fixed time step $h$ and as a basic integrator for the RTS-RK method (with uniformly distributed time steps and $p = 2$), thus computing the posterior distributions $\pi^h(y_0 \mid y_{\mathrm{obs}})$ and $\pi^h_{\mathrm{prob}}(y_0 \mid y_{\mathrm{obs}})$ defined in \cref{sec:BayesianInference}, respectively. Moreover, we compute the posterior distribution given by a non-symplectic method, the Heun's scheme, which is a classical second order method. The time step $h$ is varied for the three methods above in order to study whether the approximate posterior concentrates towards the true posterior distribution $\pi(y_0 \mid y_{\mathrm{obs}})$.

We can observe in \cref{fig:Bayes} that the posterior distributions given by Heun's method are concentrated away from the true value of the initial condition for the larger values of the time step. In fact, Heun's method is not symplectic, and a deviation on the energy $Q$ is produced when integrating the dynamical system forward in time. Hence, initial conditions with a different energy level with respect to the observation are mapped by the approximate forward model to points which are close to the observations, and as a result the posterior distribution is concentrated far from the true value. This behaviour is corrected using the Störmer--Verlet method due to its symplecticity. However, we remark that the posterior distribution for $h = 0.2$ is still concentrated on a biased value of the initial condition, without any indication of this bias given by the posterior's variance. Applying the RTS-RK method together with PMMH instead gives nested posterior distributions whose variance quantifies the uncertainty of the numerical solver. This favourable behaviour is possible due to the numerical error quantification of probabilistic methods, which has been already shown in \cite{CGS17, COS17}, together with the good energy conservation properties of the RTS-RK method when a symplectic integrator is used as its deterministic component as proved in \cref{thm:RTSHamiltonian}.

\section{Conclusion}

In this work we introduced the RTS-RK method, a novel probabilistic integrators for ODEs built on Runge--Kutta numerical integrators {with random time steps}. In particular, we analysed its weak and mean-square convergence properties, as well as the quality of Monte Carlo estimators drawn from the probabilistic solution. Geometric properties such as the conservation of first integrals and the approximation of Hamiltonians over long time intervals have been extensively treated theoretically. Finally, we showed heuristically the advantageous properties of the probabilistic approach in Bayesian inference problems with respect to the classic deterministic approach when the discretization is not in the asymptotic regime $h \to 0$. The validity of our theoretical contributions is strengthened by an extensive series of numerical examples.

The RTS-RK method is partially inspired by the probabilistic method based on additive noise perturbations presented in \cite{CGS17} and further analysed in \cite{LSS19b}, with which it shares convergence properties and the advantageous behaviour in inverse problems. Nonetheless, our method fills the void of geometry-aware probabilistic integrators for ODEs, and thus it represents a step forward in the field of probabilistic numerics for differential equations.


\section*{Appendix}


\subsection*{A modified stochastic differential equation}
In \cref{rem:sde}, we claim the existence of a modified stochastic differential equation (SDE) whose solution is well approximated by the RTS-RK method. Let us denote by $\tilde f$ the function defining the modified equation corresponding to the numerical flow $\Psi_h$ truncated after $l$ terms, i.e., 
\begin{equation}
	\tilde f(y) = f(y) + h^q f_{q+1}(y) + h^{q+1} f_{q+2}(y) + \ldots + h^l f_{l+1}(y).
\end{equation}
Details about the construction of such a function can be found in \cref{sec:Hamiltonian_2}. In particular, analyticity of the function $f$ is needed for a rigorous backward error analysis to hold. Therefore, we will refer in this section to \cref{as:RegHamiltonian} (see \cref{sec:Hamiltonian_2}). For the additive noise method presented in \cite{CGS17}, the authors consider the SDE 
\begin{equation}\label{eq:sde_AN}
	\dd Y = \tilde f(Y) \dd t + \sqrt{Q h^{2p}} \dd W,
\end{equation}
where $W$ is a $d$-dimensional standard Brownian motion. It is possible to show \cite[Theorem 2.4]{CGS17} that the solution of \eqref{eq:sde_AN} satisfies
\begin{equation}
	\abs{\E\big(\Phi(Y_N) - \Phi(Y(T)) \mid Y_0 = y \big)} \leq Ch^{2p},		
\end{equation}
where $T = Nh$ and $Y_N$ is the numerical solution given by the additive noise method after $N$ steps. Here, we present a similar construction for the RTS-RK method. In particular, let us consider the modified SDE
\begin{equation}\label{eq:sde_RTS}
	\dd\tilde Y = \Big(\tilde f(\tilde Y) + \frac{1}{2}Ch^{2p}\partial_{tt}\Psi_h(\tilde Y)\Big) \dd t + \sqrt{Ch^{2p}\partial_t \Psi_h(\tilde Y)\partial_t \Psi_h(\tilde Y)^\top} \dd W,
\end{equation}
where $C$ is given in \cref{as:hStrong}.\ref{as:hStrong_Var}. Let us denote by $\tilde \diffL$ the generator of \eqref{eq:sde_RTS}, which can be written explicitly as
\begin{equation}
	\tilde \diffL = \Big(\tilde f + \frac{1}{2}Ch^{2p}\partial_{tt}\Psi_h\Big) \cdot \nabla + \frac{1}{2}Ch^{2p}\partial_t \Psi_h \partial_t \Psi_h^\top : \nabla^2,
\end{equation}
and, adopting the semi-group notation, it satisfies
\begin{equation}
	\E \big(\Phi(\tilde Y(h))\mid \tilde Y(0) = y\big) = e^{h\tilde \diffL}\Phi(y).
\end{equation}
In the following lemma, we consider the error over one step between the numerical solution given by the RTS-RK method and the solution of \eqref{eq:sde_RTS} in the weak sense. The proof is inspired by the calculations presented in \cite[Section 2.4]{CGS17}.
\begin{lemma}\label{lem:WeakSDELocal} Under the assumptions of \cref{lem:WeakLocalOrder} and if \cref{as:RegHamiltonian} holds, then
	\begin{equation}
		\abs{\E\big(\Phi(Y_1)- \Phi(\tilde Y(h))\mid Y_0 = y\big)} \leq C h^{2p+1},
	\end{equation}
	where $C$ is {a} positive constant independent of $h$ and of $y$, $\tilde Y$ is the solution of \eqref{eq:sde_RTS} and $Y_1$ is the numerical solution given by the RTS-RK method after one step.
\end{lemma}
\begin{proof} Let us consider the modified ODE 
	\begin{equation}\label{eq:ModifiedEquation}
		\hat y'(t) = \tilde f(\hat y), 
	\end{equation}
	and denote its flow as $\hat \phi_t$. The generator $\hat \diffL = \tilde f \cdot \nabla$ satisfies, adopting the semi-group notation,
	\begin{equation}
		\Phi(\hat \phi_h(y)) = e^{h\hat \diffL}\Phi(y).
	\end{equation} 
	We can now compute the distance between the solution to \eqref{eq:sde_RTS} and \eqref{eq:ModifiedEquation} as
	\begin{equation}
	\begin{aligned}
		e^{h \tilde \diffL}\Phi(y)- e^{h \hat \diffL}\Phi(y) &= e^{h\tilde f \cdot \nabla}\Big(e^{\frac{1}{2}Ch^{{2p+1}}\partial_{tt}\Psi_h \cdot \nabla + \frac{1}{2}Ch^{{2p+1}}\partial_t \Psi_h \partial_t \Psi_h^\top : \nabla^2} - I\Big) \Phi(y)\\
		&= \big(1 + \OO(h)\big)\Big(\frac{1}{2}Ch^{{2p+1}}\partial_{tt}\Psi_h \cdot \nabla + \frac{1}{2}Ch^{{2p+1}}\partial_t \Psi_h \partial_t \Psi_h^\top : \nabla^2 + \OO\big(h^{{4p+1}}\big)\Big)\Phi(y) \\
		&= \frac{1}{2}Ch^{2p+1}\partial_{tt}\Psi_h \cdot \nabla\Phi(y) + \frac{1}{2}Ch^{2p+1}\partial_t \Psi_h \partial_t \Psi_h^\top : \nabla^2 \Phi(y) + \OO\big(h^{{4p+1}}\big).
	\end{aligned}
	\end{equation}
	Let us recall that equation \eqref{eq:DistanceProbDet} gives
	\begin{equation}
	\begin{aligned}
		e^{h\diffL_h}\Phi(y) - \Phi(\Psi_h(y)) &= \frac{1}{2} Ch^{2p+1}\partial_{tt}\Psi_h(y) \cdot \nabla\Phi(y)\\
		&\quad +\frac{1}{2}Ch^{2p+1}\partial_t \Psi_h(y) \partial_t \Psi_h(y)^\top  \colon \nabla^2\Phi(y) + \OO(h^{2p+1}),
	\end{aligned}
	\end{equation}
	which implies that
	\begin{equation}
		e^{h\tilde \diffL}\Phi(y)- e^{h\diffL_h}\Phi(y) = e^{h\hat \diffL}\Phi(y) - \Phi(\Psi_h(y)) + \OO(h^{2p+1}).
	\end{equation}	
	Now, the theory of backward error analysis (see \cref{sec:Hamiltonian_2} or e.g. \cite[Chapter IX]{HLW06}) guarantees that
	\begin{equation}
		e^{h\hat \diffL}\Phi(y) - \Phi(\Psi_h(y)) = \OO(h^{q+l+2}).
	\end{equation}
	Choosing $l = 2p - q - 1$, we have therefore
	\begin{equation}
		e^{h\tilde \diffL}\Phi(y)- e^{h\diffL_h}\Phi(y) = \OO(h^{2p+1}),
	\end{equation}
	which is the desired result.
\end{proof}
The error can be then propagated to final time as in \cref{thm:weakOrder}, as presented in the following theorem. 
\begin{theorem} Under the assumptions of \cref{lem:WeakSDELocal} and \cref{thm:weakOrder}, and if there exists a constant $L > 0$ independent of $h$ such that for all $\Phi \in \mathcal C^\infty_b(\R^d, \R)$ 
	\begin{equation}
		\sup_{u\in\R^d} \abs{e^{h\tilde \diffL}\Phi(u)} \leq (1 + Lh)\sup_{u\in\R^d}\abs{\Phi(u)},
	\end{equation}
	then it holds
	\begin{equation}
		\abs{\E\big(\Phi(Y_N)- \Phi(\tilde Y(T))\mid Y_0 = y\big)} \leq Ch^{2p},
	\end{equation}
	where $T = Nh$ and $C$ is a positive constant independent of $h$ and of $y$, $\tilde Y$ is the solution of \eqref{eq:sde_RTS} and $Y_N$ is the numerical solution given by the RTS-RK method after $N$ steps.
\end{theorem} 
\begin{proof} {The proof follows by replacing $\mathcal{L}$ with $\tilde{\mathcal{L}}$ and \cref{lem:WeakLocalOrder} with \cref{lem:WeakSDELocal} in the proof of \cref{thm:weakOrder}.}
\end{proof}

\subsection*{Proof of \cref{lem:BoundEtaJ}}
In the following, we denote by $\llbracket a, b \rrbracket$ the interval $\llbracket a, b \rrbracket = [a, b]$ if $a < b$ and $\llbracket a, b \rrbracket = [b, a]$ if $a \geq b$. Let us first consider $r \geq 2$ and the function $\gamma_r(x) = x^r e^{-r\kappa/x}$, whose first derivative is given by
	\begin{equation}
	\gamma_r'(x) = rx^{r-2}(x + \kappa) e^{-r\kappa/x}.
	\end{equation}
	Under \cref{as:RTSMoments} we have that $H_j \leq Mh$ almost surely, and hence for any $t \in \llbracket h, H_j\rrbracket$
	\begin{equation}
	\abs{\gamma_r'(t)} \leq r (Mh)^{r-2} (Mh + \kappa)e^{-r\kappa/(Mh)},
	\end{equation}
	{where we exploited that $e^{-r\kappa/x}$ is a growing function of $x$.} The fundamental theorem of calculus gives
	\begin{equation}
	\begin{aligned}
	\abs{\gamma_r(H_j)} &= \Big\lvert \gamma_r(h) + \int_{h}^{H_j} \gamma_r'(t) \dd t \,\Big\rvert\\
	&\leq \gamma_r(h) + r (Mh)^{r-2} (Mh + \kappa)e^{-r\kappa/(Mh)} \abs{H_j - h}, \quad \text{almost surely}.
	\end{aligned}
	\end{equation}
	Taking expectation on both sides and {since by \eqref{eq:BoundEtaJ} it holds} $\abs{\eta_j}^r \leq C\gamma_r(H_j)$ we obtain
	\begin{equation}
	\E\abs{\eta_j}^r \leq C\big(\gamma_r(h) + r M^{r-2} h^{p+r-3/2} (Mh + \kappa)e^{-r\kappa/(Mh)}\big),
	\end{equation} 
	which proves the desired inequality. This is because \cref{as:RTSMoments} {and \cref{as:hStrong}.\ref{as:hStrong_E} imply} that $M \geq 1$, and because $Mh$ can be bounded by {$M$}. Let us now consider $r = 1$. In this case we have for $t \in \llbracket h, H_j\rrbracket$
	\begin{equation}
	\abs{\gamma_1'(t)} \leq (mh)^{-1} (Mh + \kappa)e^{-\kappa/(Mh)}, \quad \text{almost surely.}
	\end{equation}
	Hence, we apply the same reasoning as above and obtain almost surely
	\begin{equation}
	\abs{\gamma_1(H_j)} \leq \gamma_1(h) + (mh)^{-1} (Mh + \kappa) e^{-\kappa/(Mh)} \abs{H_j - h},
	\end{equation}
	which implies the desired result by proceeding as above.  \qed

\subsection*{Proof of \cref{lem:Expansion}} 
We first expand the square as
	\begin{equation}\label{eq:ApdProofFirstExp}
	\begin{aligned}
		\Big(\sum_{j=0}^{n-1} \Big(\sum_{k=q}^{N-1} a_{jk} + b_j\Big)\Big)^2 &= \sum_{j=0}^{n-1} \Big(\sum_{k=q}^{N-1} a_{jk} + b_j\Big)^2 \\
		&+ 2 \sum_{j=1}^{n-1} \sum_{i=0}^{j-1}\Big(\sum_{k=q}^{N-1} a_{jk} + b_j\Big)\Big(\sum_{k=q}^{N-1} a_{ik} + b_i\Big).
	\end{aligned}
	\end{equation}
	Then, we expand the square in the first sum and obtain
	\begin{equation}\label{eq:ApdProofSndExp}
	\begin{aligned}
		\Big(\sum_{k=q}^{N-1} a_{jk} + b_j\Big)^2 &= \Big(\sum_{k=q}^{N-1} a_{jk}\Big)^2 + b_j^2 + 2 b_j \sum_{k=q}^{N-1} a_{jk}\\
		&= \sum_{k=q}^{N-1} a_{jk}^2 + 2\sum_{k=q+1}^{N-1}\sum_{l=q}^{k-1} a_{jk} a_{jl} + b_j^2 + 2 b_j \sum_{k=q}^{N-1} a_{jk}\\
		&= a_{jq}^2 + \sum_{k=q+1}^{N-1} a_{jk}^2 + 2\sum_{k=q+1}^{N-1}\sum_{l=q}^{k-1} a_{jk} a_{jl} + b_j^2 + 2 b_j \sum_{k=q}^{N-1} a_{jk}.
	\end{aligned}
	\end{equation}	
	We then rewrite the term appearing in the double sum in \eqref{eq:ApdProofFirstExp} as
	\begin{equation}\label{eq:ApdProofThirdExp}
	\begin{aligned}
		\Big(\sum_{k=q}^{N-1} a_{jk} + b_j\Big)\Big(\sum_{k=q}^{N-1} a_{ik} + b_i\Big) &= a_{jq}a_{iq} + \sum_{k=q}^{N-1}\sum_{\substack{l=q \\ l+k > 2q}}^{N-1} a_{jk}a_{il} \\
		&+ b_j \sum_{k=q}^{N-1} a_{ik} + b_i \sum_{k=q}^{N-1} a_{jk} + b_i b_j 
	\end{aligned}
	\end{equation}
	Substituting the expressions {\eqref{eq:ApdProofSndExp} and \eqref{eq:ApdProofThirdExp} in \eqref{eq:ApdProofFirstExp}}, we finally get 
	\begin{equation}
		\Big(\sum_{j=0}^{n-1} \Big(\sum_{k=q}^{N-1} a_{jk} + b_j\Big)\Big)^2 = \sum_{j=0}^{n-1} a_{jq}^2 + 2 \sum_{j=1}^{n-1} \sum_{i=0}^{j-1} a_{jq}a_{iq} + R(a) + S(a, b),
	\end{equation}
	where the remainder $R(a)$ can be written as $R = R_1 + R_2 + R_3$ where
	\begin{equation}\label{eq:ApdProofRemainderR}
	\begin{alignedat}{2}
		&R_1(a) = \sum_{j=0}^{n-1}\sum_{k=q+1}^{N-1} a_{jk}^2,	&&R_2(a) = 2 \sum_{j=0}^{n-1}\sum_{k=q+1}^{N-1}\sum_{l=q}^{k-1} a_{jk} a_{jl},\\
		&R_3(a) = 2 \sum_{j=1}^{n-1} \sum_{i=0}^{j-1} \sum_{k=q}^{N-1}\sum_{\substack{l=q \\ l+k > 2q}}^{N-1} a_{jk}a_{il},
	\end{alignedat}
	\end{equation}
	and the remainder $S(a, b)$ can be written as $S = S_1 + S_2 + S_3 + S_4$ where
	\begin{equation}
	\begin{alignedat}{2}
		&S_1(a, b) = \sum_{j=0}^{n-1} b_j^2,  &&S_2(a, b) = 2\sum_{j=1}^{n-1} \sum_{i=0}^{j-1} b_i b_j, \\
		&S_3(a, b) = 2\sum_{j=1}^{n-1} \sum_{k=q}^{N-1} b_j a_{jk}, \qquad &&S_4(a, b) = 2\sum_{j=1}^{n-1}\sum_{i=0}^{n-1}\Big(b_j \sum_{k=q}^{N-1} a_{ik} + b_i \sum_{k=q}^{N-1} a_{jk}\Big).\\
	\end{alignedat}
	\end{equation}
	which proves the desired result. \qed 

\subsection*{Proof of \cref{lem:Remainder}}
In the following, all the constants are independent of $h$ and $n$, but can depend on $N$ and $q$. Moreover, since $h < 1$, we often apply $h^r \leq h^s$ for $r \geq s$. We first notice that, under \cref{as:PsiStrong} and \cref{as:RegHamiltonian}, we get for all $j = 0, \ldots, n-1$ and $k = q, \ldots, N-1$
	\begin{equation}\label{eq:RemainderDeltaBound}
	\begin{aligned}
		\abs{\Delta_{j,k}} &= \abs{Q_{k+1}(Y_j) - Q_{k+1}(Y_{j+1})}\\
		&\leq C {\norm{\Psi_0(Y_j) - \Psi_{H_j}(Y_j)}}\\
		&\leq C_\Delta \abs{H_j},
	\end{aligned}
	\end{equation}
	almost surely and where $C_\Delta$ is independent of $h$. {Above, we exploited that $Q_{k+1}$ is Lipschitz continuous for all $k=q, \ldots, N+1$ due to \cref{as:RegHamiltonian}.} Let us now consider $R(\Delta)$. {Due to \eqref{eq:RemainderDeltaBound} and to \cref{as:RTSMoments}}, we have
	\begin{equation}\label{eq:RemainderFirstBound}
	\begin{aligned}
		\E (H_j^k - h^k)^2 \Delta_{j,k}^2 &\leq {C_\Delta^2}\E(H_j^{k+1} - H_j h^k)^2\\
		&= {C_\Delta^2}\big(h^{2(k+1)} + C_{2(k+1)}h^{2p+2(k+1)-1} + h^{2(k+1)} + C_2 h^{2p+2k+1}  \\
		&\quad - 2h^{2k+2} - 2C_{k+2}h^{2p+2k+1}\big)\\
		&= {C_\Delta^2}\big((C_{2(k+1)} + C_2 - 2C_{k+2})h^{2p+2k+1}\big)\\
		&\leq C h^{2p+2k+1},
	\end{aligned}
	\end{equation}
	where $C > 0$ is a positive constant. Now, since $k \geq q+1$, we get
	\begin{equation}\label{eq:RemainderSecondBound}
		\E (H_j^k - h^k)^2 \Delta_{j,k}^2 \leq C h^{2(p+q+1)}.
	\end{equation}
	Hence, for $R_1(\Delta)$ there exists a constant $\tilde C_1$ such that
	\begin{equation}
		\E R_1(\Delta) \leq \tilde C_1 n h^{2(p+q+1)}.
	\end{equation}
	We now proceed to the second remainder $R_2(\Delta)$. Applying the Cauchy--Schwarz inequality and \eqref{eq:RemainderFirstBound} we get
	\begin{equation}
	\begin{aligned}
		 \E \big((H_j^k - h^k)\Delta_{j,k}(H_j^l - h^l) \Delta_{j,l}\big) &\leq \Big(\E \big((H_j^k - h^k)^2\Delta_{j,k}^2\Big)^{1/2} \Big(\E \big((H_j^l - h^l)^2\Delta_{j,l}^2\Big)^{1/2}\\
		 &\leq C h^{2p+k+l+1},
	\end{aligned}
	\end{equation}
	where $C > 0$ is a positive constant. Now, since {in the definition of $R_2(a)$ in \eqref{eq:RemainderDeltaBound} we have $k \geq q+1$ and $l \geq q$, we have here $k+l \geq 2q+1$}. Therefore, there exists a constant $\tilde C_2$ such that
	\begin{equation}
		\E R_2(\Delta) \leq \tilde C_2 n h^{2(p+q+1)}.
	\end{equation}	
	We now consider the term $R_3(\Delta)$. Since $H_i$ and $H_j$ are independent for $i \neq j$, we have
	\begin{equation}
	\begin{aligned}
		\E \big((H_j^k - h^k)\Delta_{j,k}(H_i^l - h^l) \Delta_{i,l}\big) &= \E (H_j^k - h^k)\Delta_{j,k} \E (H_i^l - h^l) \Delta_{i,l}.
	\end{aligned}
	\end{equation}
	Computing the two factors singularly, we have {due to \eqref{eq:RemainderDeltaBound} and to \cref{as:RTSMoments}}
	\begin{equation}\label{eq:RemainderThirdBound}
	\begin{aligned}
		\E (H_j^k - h^k)\Delta_{j,k} &\leq C_\Delta \E(H_j^{k+1} - H_jh^k) \\
		&= C_\Delta C_{k+1} h^{2p+k},
	\end{aligned}
	\end{equation}
	and analogously for $\E (H_i^l - h^l) \Delta_{i,l}$. Then, since $k+l \geq 2q+1$
	\begin{equation}\label{eq:RemainderFourthBound}
		\E \big((H_j^k - h^k)\Delta_{j,k}(H_i^l - h^l) \Delta_{i,l}\big) \leq C_\Delta^2 C_{k+1} C_{l+1} h^{2(2p+q+1/2)}.
	\end{equation}
	Hence, we have for a constant $\tilde C_3 > 0$ 
	\begin{equation}
		\E R_3(\Delta) \leq \tilde C_3 n^2 h^{2(2p+q+1/2)}.
	\end{equation}
	Finally, replacing $t_n = nh$, we can write for a constant $C > 0$
	\begin{equation}
	\begin{aligned}
		\E R(\Delta) &\leq (\tilde C_1 + \tilde C_2) n h^{2(p+q+1)} + \tilde C_3 n^2 h^{2(2p+q+1/2)}\\
		&= (\tilde C_1 + \tilde C_2) t_n h^{2(p + q + 1/2)} + \tilde C_3 t_n^2 h^{2(2p + q - 1/2)}.
	\end{aligned}
	\end{equation}
	Let us now consider $S(\Delta, \eta)$. First, we notice that under the assumption $p \geq 3/2$ we have for any $r \geq 1$, $\min\{r, p+r-3/2\} = r$, and therefore \cref{lem:BoundEtaJ} simplifies to 
	\begin{equation}
		\E\abs{\eta_j}^r \leq Ch^r e^{-r\kappa/(Mh)}.
	\end{equation}
	We first consider $S_1(\Delta, \eta)$. Applying \cref{lem:BoundEtaJ} with $r = 2$, we obtain for a constant $\hat C_1 > 0$ 
	\begin{equation}
		\E S_1(\Delta, \eta) \leq \hat C_1 n h^2 e^{-2\kappa/(Mh)}.
	\end{equation}
	For the second term $S_2(\Delta, \eta)$, we have {by \eqref{eq:BoundEtaJ}} that $\abs{\eta_i}\leq CH^i e^{-\kappa/H_i}$ and $\eta_j\leq CH^j e^{-\kappa/H_j}$ almost surely. These two bounds are independent for $i \neq j$ and therefore, applying \cref{lem:BoundEtaJ} with $r = 1$, we have for a constant $\hat C_2 > 0$
	\begin{equation}
		\E S_2(\Delta, \eta) \leq \hat C_2 n^2 h^2 e^{-2\kappa/(Mh)}.
	\end{equation}
	We now consider the third remainder $S_3(\Delta, \eta)$. Applying the Cauchy--Schwarz inequality, we obtain
	\begin{equation}
		\E \eta_j (H_j^k - h^k) \Delta_{j,k} \leq (\E\eta_j^2)^{1/2} (\E(H_j^k - h^k)^2 \Delta_{j,k}^2)^{1/2}.
	\end{equation}
	Applying \cref{lem:BoundEtaJ} with $r = 2$ to the first factor and \eqref{eq:RemainderFirstBound} to the second we get
	\begin{equation}
	\begin{aligned}
		\E \eta_j (H_j^k - h^k) \Delta_{j,k} &\leq C h e^{-\kappa/(Mh)}h^{p+k+1/2}\\
		&= C h^{p+k+3/2}e^{-\kappa/(Mh)}
	\end{aligned}
	\end{equation}
	Now, since $k \geq q$, we have for a constant $\hat C_3 > 0$
	\begin{equation}
		\E S_3(\Delta, \eta) \leq \hat C_3 n h^{p+q+3/2}e^{-\kappa/(Mh)}. 
	\end{equation}
	Finally, we consider the last term $S_4(\Delta, \eta)$. Since {by \eqref{eq:BoundEtaJ} it holds $\abs{\eta_j} \leq CH_j e^{-\kappa/H_j}$} almost surely, and this bound is independent of $H_i$ for $i \neq j$, applying \eqref{eq:RemainderThirdBound} and \cref{lem:BoundEtaJ} we have
	\begin{equation}
	\begin{aligned}
		\E\eta_j(H_i^k - h^k)\Delta_{i,k} &= \E\eta_j \E (H_i^k - h^k)\Delta_{i,k}\\
		&\leq Ch e^{-\kappa/(Mh)}h^{2p+k},
	\end{aligned}
	\end{equation}
	which, since $k \geq q$, implies that there exists a constant $\hat C_4 > 0$ such that 
	\begin{equation}
		\E S_4(\Delta, \eta) \leq \hat C_4 n^2 h^{2p+q+1}e^{-\kappa/(Mh)}.
	\end{equation}
	Finally, replacing $t_n = nh$, we can write
	\begin{equation}
	\begin{aligned}
		\E S(\Delta, \eta) &\leq (\hat C_1 n h^2 + \hat C_2 n^2 h^2) e^{-2\kappa/(Mh)} + (\hat C_3 n h^{p+q+3/2} + \hat C_4 n^2 h^{2p+q+1})e^{-\kappa/(Mh)}\\
		&= (\hat C_1 t_n h + \hat C_2 t_n^2) e^{-2\kappa/(Mh)} + (\hat C_3 t_n h^{p+q+1/2} + \hat C_4 t_n^2 h^{2p+q-1})e^{-\kappa/(Mh)},
	\end{aligned}
	\end{equation}
	which completes the proof.\qed
	
\def\cprime{$'$}

\end{document}